\documentclass[12pt,leqno]{article}
\usepackage{amssymb}

\newtheorem{theorem}{Theorem}[section]
\newtheorem{lemma}[theorem]{Lemma}
\newtheorem{remark}[theorem]{Remark}
\newtheorem{example}[theorem]{Example}
\newtheorem{corollary}[theorem]{Corollary}

\def\square{{\vcenter{\vbox{\hrule height.3pt
                  \hbox{\vrule width.3pt height5pt \kern5pt
                     \vrule width.3pt}
                  \hrule height.3pt}}}}
\def\qed{{\hfill $\square$ \smallskip}}

\def \Proof{\noindent{\bf Proof\mbox{\quad}}}

\def \be{\begin{equation}}
\def \ee{\end{equation}}
\def \bt{\begin{theorem}}
\def \et{\end{theorem}}
\def \bc{\begin{corollary}}
\def \ec{\end{corollary}}
\def \bex{\begin{example}}
\def \eex{\end{example}}

\def \br{\begin{remark}}
\def \er{\end{remark}}
\def \bl{\begin{lemma}}
\def \el{\end{lemma}}

\def \bea{\begin{eqnarray}}
\def \eea{\end{eqnarray}}
\def \bas{\begin{eqnarray*}}
\def \eas{\end{eqnarray*}}

\def \ga{\gamma}
\def \Ga{\Gamma}
\def \de{\delta}

\def \ep{\epsilon}

\def \la{\lambda}
\def \La{\Lambda}

\def \ff{\infty}
\def \wh{\widehat}
\def \wt{\widetilde}

\def \rar{\rightarrow}

\def \R{{\bf R}}

\def \Z{{\bf Z}}

\def \({\left(}
\def \){\right)}

\def \lc{\left\{}
\def \rc{\right\}}

\def \nn{\nonumber}

\def \bs{\begin{slide} }
\def \es{\end{slide} }

\def \bpr{\begin{proof} }
\def \epr{\end{proof} }

\def \br{\begin{remark}}
\def \er{\end{remark}}
\def \bl{\begin{lemma}}
\def \el{\end{lemma}}

\newcommand{\eqnsection}{
         \renewcommand{\theequation}{\thesection.\arabic{equation}}
         \makeatletter
         \csname @addtoreset\endcsname{equation}{section}
         \makeatother}

\eqnsection

\def \R{{\Bbb R}}
\def \E{{{\Bbb E}\,}}
\def \P{{\Bbb P}}
\def \Q{{\Bbb Q}}
\def \Z{{\Bbb Z}}
\def\F{{\cal F}}

\def\lam{{\lambda}}

\def\proof{{\medskip\noindent {\bf Proof. }}}

\def\eps{\varepsilon}

 \def\qq {\qquad}
\def\frac#1#2{{#1\over #2}}

\def\wt{\widetilde}
\def\ol{\overline}
\def\wh{\widehat}

\def\ni{\noindent }

\def\bs{\bigskip}

\def\tfrac#1#2{{\textstyle {#1\over #2}}}

\def\sA {{\cal A}}

\def\sS {{\cal S}}

\begin{document}

\title{Moderate deviations and laws of the iterated logarithm for the
renormalized self-intersection local times of planar random walks}
\author{ Richard F.~Bass \thanks {Research partially supported by NSF
grant DMS-0244737.}\hspace{ .2in}  Xia Chen \thanks { Research partially
supported by NSF grant DMS-0405188.}\hspace{ .2in}  Jay Rosen\thanks
      {This research was  supported, in part, by grants from the National
Science Foundation and PSC-CUNY.}}

\maketitle

\begin{abstract} We study moderate deviations for the renormalized
self-intersection local time of planar random walks. We also prove laws of
the iterated logarithm for such local times.
\end{abstract}

\section{ Introduction}

      Let $\{S_n\}$ be a symmetric random walk on $\Z^2$
      with covariance matrix $\Gamma$.  Let
\medskip
\begin{equation}
      B_n=\sum_{1\le j<k\le n}\de( S_{ j},S_{ k})\label{1.1}
\end{equation} where
\begin{equation}
\de( x,y)=\left\{\begin{array}{ll} 1&\mbox{  if $x=y$}\\ 0&\mbox{
otherwise}
\end{array}\right.\label{1.0}
\end{equation} is the usual Kroenecker delta. We refer to $B_n$ as the
self-intersection local time up to time $n$. We call
\[ \ga_{ n}=:B_n-\E B_n
\] the renormalized self-intersection local time of the random walk up to
time
$n$.

In
\cite{BRilt} it was shown that $\ga_{ n}$, appropriately scaled, converges
to the renormalized self-intersection local time of planar Brownian motion.
Renormalized self-intersection local time for Brownian motion was
originally studied by Varadhan
\cite{Va} for its role in quantum field theory. Renormalized
self-intersection local time turns out to be the right  tool  for the solution
of certain ``classical'' problems such as the asymptotic expansion of the
area of the Wiener sausage in the plane and the range of random walks,
\cite{BR}, \cite{LGa}, \cite{LG}.

One of the applications of self-intersection local time is to  polymer
growth. If $S_n$ is a planar random walk and
$\P$ is its law, one can construct self-repelling and self-attracting random
walks by defining
\[ d\Q_n/d\P= c_n e^{\zeta B_n/n}, \] where $\zeta$ is a parameter and
$c_n$ is chosen to make $\Q_n$ a probability measure. When
$\zeta<0$, more weight is given to those paths with a small number of
self-intersections, hence $\Q_n$ is a model for a self-repelling random
walk. When $\zeta>0$, more weight is given to paths with a large number
of self-intersections, leading to a self-attracting random walk. Since $\E
B_n$ is deterministic, by modifying  $c_n$, we can write
\[ d\Q_n/d\P=c_n e^{\zeta(B_n-\E B_n)/n}.
\]
It is known that for small positive $\zeta$ the self-attracting
random walk grows with $n$ while for large $\zeta$ it ``collapses,''
and its diameter remains bounded in mean square. It has been an
open problem to determine the critical value of $\zeta$ at which the
phase transition takes place. The work \cite{BC} suggested that the
critical value $\zeta_c$ could be expressed in terms of the best
constant of a certain Gagliardo-Nirenberg inequality, but that work
was for planar Brownian motion, not for random walks. In the current
paper we obtain moderate deviations estimates for $\gamma_n$ and these
are in terms of the best constant of the Gagliard-Nirenberg inequality; see
Theorem \ref{theo-md}.  However the critical constant $\zeta_c$ is
different (see Remark \ref{largedev}) and it is still an open problem
to determine it.
See \cite{BoS} and \cite{BS} for details and
further information on these models.

In the present paper we study moderate deviations of $\ga_{ n}$. Before
stating our main theorem we  recall one of the Gagliardo-Nirenberg
inequalities:
\[ \|f\|_{ 4}\leq C  \|\nabla f\|_{ 2}^{ 1/2} \|f\|_{ 2}^{ 1/2},
\] which is valid for $f\in C^1$ with compact support, and can then be
extended to more general $f$'s. We define
$\kappa (2,2)$ to be the infimum of those values of $C$  for which the
above inequality holds. In particular,
$0<\kappa (2,2)<\ff$. For further details, see \cite{C04}.

In this paper we will always assume that the smallest group which supports
$\{S_n\}$ is
$\Z^2$.  For simplicity we assume  further that our random walk is
strongly aperiodic.

\bt\label{theo-md}  Let $\{b_n\}$ be a positive sequence satisfying
\be\label{restrict}
\lim_{ n\rar
\ff}b_n=\infty\hskip.1in\hbox{and}\hskip.1in b_n=o(n).\label{18.17}
\ee
For any $\lambda>0$,
\begin{equation}
\lim_{n\to\infty}{1\over b_n}\log\P\Big\{B_n-\E B_n\ge \lambda
nb_n\Big\} =-\lambda\sqrt{\det \Gamma}\,\,\kappa (2,2)^{-4}.\label{1.2}
\end{equation}
\et

We call Theorem \ref{theo-md} a moderate deviations theorem rather than
a large deviations result because of the  second restriction in
(\ref{restrict}). Our techniques do not apply when this restriction is
not present, and
and in fact it is not hard to show that the value on the right hand side of
(\ref{1.2}) should be different when $b_n\approx n$; see Remark \ref{largedev}.

Moderate deviations for $-\ga_{ n}$ are more subtle. In the next theorem
we obtain the correct rate, but not the precise constant.

\bt\label{theo-ld} Suppose $\E|S_{ 1}|^{ 2+\de}<\ff$ for some
$\de>0$. There exist
$C_1, C_2>0$ such that for any
$\theta >0$ and sequence
$ b_n\to\infty\hskip.1in\hbox{with}\hskip.1in b_n =o(n^{1/\theta})$
\begin{eqnarray} -C_1&\leq&
\liminf_{n\to\infty}b_n^{-\theta}\log\P\Big\{\E B_n -B_n
\ge\theta (2\pi)^{-1}\det(\Gamma)^{-1/2}n\log b_n\Big\}
\nn\\ &\leq& \limsup_{n\to\infty}b_n^{-\theta}\log\P\Big\{\E B_n -B_n
\ge\theta (2\pi)^{-1}\det(\Gamma)^{-1/2}n\log b_n\Big\}\nn\\ &\le&
-C_2.
\label{1.3}
\end{eqnarray}
\et

Here are the corresponding  laws of the iterated logarithm for $\ga_{ n}$.

\bt\label{theo-lil}
\begin{equation}
\limsup_{n\to\infty}{B_n -\E B_n\over n\log\log
n}=\det(\Gamma)^{-1/2}
\kappa(2,2)^4\hskip.2in a.s.\label{1.4}
\end{equation} and if  $\E|S_{ 1}|^{ 2+\de}<\ff$ for some $\de>0$,
\begin{equation}
\liminf_{n\to\infty}{B_n -\E B_n\over n\log\log\log n}=-(2\pi)^{-1}
\det(\Gamma)^{-1/2}\hskip.2in a.s.\label{1.5}
\end{equation}
\et

In this paper we  deal exclusively with the case where the dimension
$d$ is 2. We note that in dimension $1$ no renormalization is needed,
which makes the results much simpler. See \cite{M, CL}. When $d\ge 3$,
the renormalized intersection local time is in the domain of attraction of a
centered normal random variable. Consequently the  tails of the weak limit
are expected to  be of Gaussian type, and  in particular, the tails  are
symmetric; see
\cite{LG}.

Theorems \ref{theo-md}-\ref{theo-lil} are the analogues of the theorems
proven in
\cite{BC} for the renormalized self-intersection local time of planar
Brownian motion. Although the  proofs for the random walk case have
some elements in common with those for  Brownian motion, the random
walk case is considerably more difficult. The major difficulty
      is the fact that we do not have Gaussian  random variables.
Consequently, the argument for the lower bound of Theorem 1.1 needs to
be very different from the one given in \cite[Lemma 3.4]{BC}. This
requires several new tools, such as Theorem
\ref{theo-exprw}, which we expect will have applications beyond the
specific needs of this paper.

\section{ Integrability}

Let $\{S'_n\}$ be an independent copy of the random walk
$\{S_n\}$. Let
\begin{equation} I_{m,n}=\sum_{ j=1}^{ m}\sum_{ k=1}^{ n}
\de( S_{ j},S'_{ k})\label{2.1r}
\end{equation} and set $I_{n}=I_{n,n}$.
      Thus
\begin{equation} I_n=\#\{(j,k)\in [1,n]^{ 2};\hskip.1in S_j=S_k'\}.
\label{2.1}
\end{equation}

\bl\label{lem-firstmom}
\begin{equation}
\E I_{m,n}\leq c\(( m+n )\log ( m+n )-m\log m-n\log n\).\label{2.0x}
\end{equation} In particular
\begin{equation}
\E(  I_n)\le c n.
\label{2.2}
\end{equation} We also have
\begin{equation}
\E I_{m,n}\leq c\sqrt{mn}.\label{2.0y}
\end{equation}
\el

\Proof Using symmetry and independence
\begin{eqnarray}
\E I_{m,n}&=&\sum_{ j=1}^{ m}\sum_{ k=1}^{ n}
\E \de( S_{ j},S'_{ k})\label{2.1s}\\ &=&\sum_{ j=1}^{ m}\sum_{ k=1}^{ n}
\E \de( S_{ j}-S'_{ k},0)\nn\\ &=&\sum_{ j=1}^{ m}\sum_{ k=1}^{ n}
\E \de( S_{ j+k},0)=\sum_{ j=1}^{ m}\sum_{ k=1}^{ n}p_{j+k }( 0).\nn
\end{eqnarray} By \cite[p. 75]{S},
\begin{equation} p_{ m}( 0)={1\over 2\pi\sqrt{\det\Gamma}}\,{ 1\over
m}\,+\, o\({ 1\over m} \)\label{4.3g1}
\end{equation} so that
\begin{equation}
\E I_{m,n}\leq c\sum_{ j=1}^{ m}\sum_{ k=1}^{ n}{1 \over j+k }
\leq c\int_{r=0}^{ m}\int_{ s=0}^{ n}{1 \over r+s }\,dr\,ds\label{2.1xx}
\end{equation} and (\ref{2.0x}) follows. (\ref{2.1}) is then immediate.
(\ref{2.0y}) follows from (\ref{2.1xx}) and the bound $(r+s )^{ -1}\leq
(\sqrt{rs})^{ -1}$.
\qed

It follows from the proof of \cite[Lemma 5.2]{C04} that for any integer
$k\geq 1$
\begin{equation}
\E(  I^{ k}_n)\le ( k!)^{ 2}(1+\E(  I_n))^{ k}\label{2.1a}
\end{equation} Furthermore, by  \cite[(5.k)]{LG} we have that
$I_n/n$ converges in distribution to a random variable with finite
moments. Hence for any  integer $k\geq 1$
\begin{equation}
\lim_{ n\rar \ff}{ \E(  I^{ k}_n)\over n^{ k}}= c_{ k}<\ff.\label{2.1b}
\end{equation}

\bl\label{lem-ubound}  There is a constant $c>0$ such that
\begin{equation}
\sup_n\E\exp\Big\{{c\over n}I_n\Big\}<\infty.
\label{2.3}
\end{equation}
\el

\proof For any
$m\ge 1$ write $l(m, n)=[n/m]+1$. Using  \cite[Theorem 5.1]{C04}
      with $p=2$ and $a=m$, and then  (\ref{2.2}), (\ref{2.1a}) and
(\ref{2.1b}), we obtain
\begin{eqnarray} \big(\E I_n^m\big)^{1/2}&\leq &
\sum_{\stackrel{ k_1+\cdots +k_m =m}{ k_1,\cdots, k_m\ge 0}}{m!\over
k_1!\cdots k_m!}
\big(\E I_{l(m,n)}^{k_1}\big)^{1/2}\cdots
\big(\E I_{l(m,n)}^{k_m}\big)^{1/2}
\nn\\ &\leq & \sum_{\stackrel{ k_1+\cdots +k_m =m}{ k_1,\cdots,
k_m\ge 0}}{C^mm!\over k_1!\cdots k_m!}k_1!\cdots k_m!
\big(\E I_{l(m,n)}\big)^{k_1/2}\cdots
\big(\E I_{l(m,n)}\big)^{k_m/2}\label{2.4}\\ &\leq &
      { 2m-1\choose m}m! C^m
\Big({n\over m}\Big)^{ m/2}
\le  { 2m\choose m}m!C^m \Big({n\over m}\Big)^{m/2}\nn
\end{eqnarray} where $C>0$ can be chosen independently of $m$ and
$n$. Hence
\begin{equation}
\E I_n^m\le  { 2m\choose m}^2C^m(m!)^2\Big({n\over m}\Big)^m
\le { 2m\choose m}^2C^m m! n^m.
\label{2.5}
\end{equation}
      Notice that
\begin{equation} { 2m\choose m}\le 4^m.
\label{2.6}
\end{equation} The conclusion then follows using the power series for
$e^{ x}$.
\qed

For any random variable $X$ we define
\[ \ol X=:X-\E X.
\] We write
\begin{equation}(m, n]^2_{ <}=\{(j,k)\in (m, n]^{ 2};\hskip.1in j<k\}
\label{4.4}
\end{equation} For any $A\subset \big\{(j,k)\in (\Z^+)^2;\hskip.1in
j<k\big\}$, write
\begin{equation}
      B(A)=\sum_{(j,k)\in A}\de(S_j,S_k )
\label{2.9}
\end{equation}  In our proofs we will use several decompositions of
$B_n$. If $J_1, \ldots, J_\ell$ are consecutive disjoint blocks of integers
whose union is $\{1, \ldots, n\}$, we have
\[ B_n=\sum_i B((J_i\times J_i)\cap (0,n]_{<})+\sum_{i<j} B(J_i\times
J_j) \] and also
\[ B_n=\sum_i B((J_i\times J_i)\cap (0,n]_{<})+\sum_{i}
B(\cup_{j=1}^{i-1} J_j)\times J_i). \]

\bl\label{lem-ibound}  There is a constant $c>0$ such that
\begin{equation}
\sup_n\E\exp\Big\{{c\over n}\vert\, \ol B_n\,\vert\Big\}<\infty.
\label{2.7}
\end{equation}
\el

\proof We first prove that there is
$c>0$ such that
\begin{equation} M\equiv \sup_n\E\exp\Big\{{c\over 2^n}\vert\, \ol
B_{2^n}\,\vert\Big\}<\infty.
\label{2.8}
\end{equation} We have
\begin{eqnarray} &&\hspace{ .2in} B_{2^n}
\label{2.9a}\\ &&
=\sum_{j=1}^n\sum_{k=1}^{2^{j-1}}B\Big(\big((2k-2)2^{n-j},\hskip.03in
(2k-1) 2^{n-j}\big]\times\big((2k-1)2^{n-j},\hskip.03in
(2k)2^{n-j}\big]\Big)\nonumber
\end{eqnarray}
      Write
\begin{eqnarray} &&\hspace{ .2in}\alpha_{j,k}=
      B\Big(\big((2k-2)2^{n-j},\hskip.03in (2k-1)
2^{n-j}\big]\times\big((2k-1)2^{n-j},\hskip.03in
(2k)2^{n-j}\big]\Big)\label{2.10}\\ && \hspace{ .4in}-\E
B\Big(\big((2k-2)2^{n-j},\hskip.03in (2k-1)
2^{n-j}\big]\times\big((2k-1)2^{n-j},\hskip.03in (2k)2^{n-j}\big]\Big)
\nonumber
\end{eqnarray} For each $1\le j\le n$, the random variables
$\alpha_{j,k}$, $\,\,k=1,\cdots, 2^{j-1}$ are i.i.d. with common
distribution
$I_{2^{n-j}}-\E I_{2^{n-j}}$. By the previous lemma there exists
$\delta >0$ such that
\begin{equation}
\sup_n\sup_{j\leq n}\E\exp\Big\{\delta {1\over 2^{n-j}}\big\vert
\alpha_{j,1}
\big\vert\Big\}<\infty.
\label{2.11}
\end{equation} By  \cite[Lemma 1]{BCR05}, there exists $\theta >0$ such
that
\begin{eqnarray} C(\theta) &\equiv &\sup_{n}\sup_{j\leq
n}\E\exp\bigg\{\theta 2^{j/2}{1\over 2^n}\Big\vert
\sum_{k=1}^{2^{j-1}}\alpha_{j,k}\Big\vert\bigg\}
\label{2.12}\\ &=& \sup_{n}\sup_{j\leq n}\E\exp\bigg\{\theta
2^{-j/2}{1\over 2^{n-j}}
\Big\vert
\sum_{k=1}^{2^{j-1}}\alpha_{j,k}
\Big\vert\bigg\}<\infty.\nonumber
\end{eqnarray} Write
\begin{equation}
\lambda_N= \prod_{j=1}^N\big(1-2^{-j/2}\big)\hskip.1in\hbox{and}
\hskip.1in \lambda_\infty = \prod_{j=1}^\infty \big(1-2^{-j/2}\big).
\label{2.13}
\end{equation}
      Using H\"older's inequality with $1/p=1-2^{-n/2},\,1/q=2^{-n/2}$ we
have
\begin{eqnarray} &&
\E\exp\Big\{\lambda_n{\theta\over 2^n}\Big\vert
\sum_{j=1}^n\sum_{k=1}^{2^{j-1}}\alpha_{j,k}\Big\vert\Big\}
\label{2.14}\\ && \le\bigg(\E\exp\Big\{\lambda_{n-1}{\theta\over
2^n}\Big\vert
\sum_{j=1}^{n-1}\sum_{k=1}^{2^{j-1}}
\alpha_{j,k}\Big\vert\Big\}\bigg)^{1-2^{-n/2}}\nn\\ &&\hspace{ 1in}
\times \bigg(\E\exp\Big\{2^{n/2}\lambda_n{\theta\over 2^n}\Big\vert
\sum_{k=1}^{2^{n-1}}\alpha_{n,k}\Big\vert\Big\}\bigg)^{2^{-n/2}}
\nonumber\\ &&\le\E\exp\Big\{\lambda_{n-1}{\theta\over
2^n}\Big\vert
\sum_{j=1}^{n-1}\sum_{k=1}^{2^{j-1}}
\alpha_{j,k}\Big\vert\Big\}C(\theta)^{2^{-n/2}}\nn
\end{eqnarray}
      Repeating this procedure,
\begin{eqnarray} &&\E\exp\Big\{\lambda_n{\theta\over 2^n}\Big\vert
\sum_{j=1}^n\sum_{k=1}^{2^{j-1}}\alpha_{j,k}\Big\vert\Big\}
\label{2.15}\\ && \le C(\theta)^{2^{-1/2}+\cdots +2^{-n/2}}
\le C(\theta)^{2^{-1/2} (1-2^{-1/2})^{-1}}\nonumber
\end{eqnarray} So we have
\begin{equation}
\sup_n\E\exp\Big\{\lambda_\infty{\theta\over 2^n}\vert \,\ol
B_{2^n}\,\vert\Big\}<\infty
\label{2.16}
\end{equation}

We now prove our lemma for general $n$. Given an integer $n\ge 2$, we
have the following unique representation:
\begin{equation}  n=2^{m_1}+2^{m_2}+\cdots +2^{m_l}
\label{2.17}
\end{equation} where $m_1>m_2>\cdots m_l\ge 0$ are integers. Write
\begin{equation}  n_0=0\hskip.1in \hbox{and}\hskip .1in
n_i=2^{m_1}+\cdots +2^{m_i},
\qq   i=1,\cdots, l.
\label{2.18}
\end{equation} Then
\begin{eqnarray} \sum_{1\le j<k\le n}\de(S_j,S_k )&=&
\sum_{i=1}^l\sum_{n_{i-1}<j< k\le n_i}\de(S_j,S_k )
+\sum_{i=1}^{l-1}B\big((n_{i-1},n_i]\times(n_i,n]\big)
\nn\\ &=& :\sum_{i=1}^lB^{ ( i)}_{ 2^{ m_{
i}}}+\sum_{i=1}^{l-1}A_i.\label{2.19}
\end{eqnarray}
      By H\"older's inequality, with $M$ as in (\ref{2.8})
\begin{eqnarray} &&\E\exp\Big\{{c\over n}\Big\vert\sum_{i=1}^l (B^{ (
i)}_{ 2^{ m_{ i}}}-\E B^{ ( i)}_{ 2^{ m_{ i}}})\Big\vert\Big\}
\label{2.20}\\ && \le \prod_{i=1}^l\bigg(\E\exp\Big\{{c\over 2^{m_i}}
\vert B^{ ( i)}_{ 2^{ m_{ i}}}-\E B^{ ( i)}_{ 2^{ m_{ i}}}\vert
\Big\}\bigg)^{2^{m_i}\over n}\le \prod_{i=1}^l M^{2^{m_i}/
n}=M.\nonumber
\end{eqnarray}

      Using  H\"older's inequality,
\begin{equation}
\E\exp\Big\{{c\over n}\sum_{i=1}^{l-1}A_i\Big\}
\le \prod_{i=1}^{l-1}\bigg(\E\exp\Big\{{c\over 2^{m_i}}A_i
\Big\}\bigg)^{2^{m_i}\over n}.
\label{2.22}
\end{equation}
      Notice that for each $1\le i\le l-1$,
\begin{equation} A_i\buildrel d\over
=\sum_{j=1}^{2^{m_i}}\sum_{k=1}^{n-n_i}
\delta(S_j,S_k')\le \sum_{j=1}^{2^{m_i}}\sum_{k=1}^{2^{m_i}}
\delta(S_j,S_k'),
\label{2.23}
\end{equation}
      where the inequality follows from
\begin{equation}  n-n_i=2^{m_{i+1}}+\cdots +2^{m_l}\le 2^{m_i}.
\label{2.24}
\end{equation} Using (\ref{2.23}) and Lemma \ref{lem-firstmom}, we can
take $c>0$ so that
\begin{equation}
\E\exp\Big\{{c\over 2^{m_i}}A_i\Big\}\le
\sup_n\E\exp\Big\{{c\over n}I_n\Big\}\equiv N<\infty.
\label{2.25}
\end{equation}
      Consequently,
\begin{equation}
\E\exp\Big\{{c\over n}\sum_{i=1}^{l-1}A_i\Big\}
\le \prod_{i=1}^{l-1}N^{2^{m_i}/n}\le N.
\label{2.26}
\end{equation} In particular, this shows that
\begin{equation}
\E\Big\{{c\over n}\sum_{i=1}^{l-1}A_i\Big\}
\le N.
\label{2.26a}
\end{equation} Combining (\ref{2.26})  and (\ref{2.26a}) with (\ref{2.20})
we have
\begin{equation}
\sup_n\E\exp\Big\{{c\over 2n}\vert \ol B_n\vert\Big\}<\infty.
\label{2.27}
\end{equation}
\qed

\bl\label{lem-expect}
\begin{equation}
\E B_n={1\over 2\pi\sqrt{\det\Gamma}}n\log n +o(n\log n),
\label{4.3}
\end{equation} and if $\E|S_{ 1}|^{ 2+2\de}<\ff$ for some $\de>0$ then
\begin{equation}
\E B_n={1\over 2\pi\sqrt{\det\Gamma}}n\log n +O(n).
\label{4.3s}
\end{equation}
\el

\proof
\begin{equation}
\E B_n=\E \sum_{1\le j<k\le n}\de( S_{ j},S_{ k}) =\sum_{1\le j<k\le
n}p_{ k-j}( 0)\label{4.3f}
\end{equation} where $p_{m}( x)=\E (  S_{ m}=x)$. If $\E|S_{ 1}|^{
2+2\de}<\ff$, then by \cite[Proposition 6.7]{Lawler},
\begin{equation} p_{ m}( 0)={1\over 2\pi\sqrt{\det\Gamma}}\,{ 1\over
m}\,+\, o\({ 1\over m^{ 1+\de}} \).\label{4.3g}
\end{equation} Since the last term is summable, it will contribute
$O(n)$ to (\ref{4.3f}). Also,
\begin{equation}
\sum_{1\le j<k\le n}{ 1\over k-j}=\sum_{ m=1}^{ n}\sum_{ i=1}^{ n-m}{
1\over m}=\sum_{ m=1}^{ n}{ n-m\over m}=n\sum_{ m=1}^{ n}{1\over
m}-n\label{4.3h}
\end{equation} and our Lemma follows from the well known fact that
\begin{equation}
\sum_{ m=1}^{ n}{1\over m}=\log n  +\ga+O\({1\over n}\)\label{4.3i}
\end{equation} where $\ga$ is Euler's constant.

If we only assume finite second moments, instead of (\ref{4.3g}) we use
(\ref{4.3g1})
      and proceed as above.
\qed

\bl\label{lem-idbound} For any $\theta >0$
\begin{equation}
\sup_n\E\exp\Big\{{\theta\over n}(\E B_n- B_n)\Big\}<\infty
\label{2.28}
\end{equation} and for any $\lambda>0$
\begin{equation}
\lim_{n\to\infty}{1\over b_n}\log\P\Big\{\E B_n -B_n
\ge \lambda nb_n\Big\} =-\infty.
\label{3.3}
\end{equation}
\el

\proof By Lemma \ref{lem-ibound} this is true for some
$\theta_o>0$. For any $\theta >\theta_o$, take an integer $m\ge 1$ such
that
$\theta m^{-1}<\theta_o$. We can write any $n$ as $n=rm+i$ with
$1\leq i<m$. Then
\begin{eqnarray} &&\hspace{ .2in}\E B_{n}- B_{n}
\label{2.30}\\ &&\le \sum_{j=1}^m\Big[\E
\sum_{(j-1)r<k,l\le jr}\de(S_k,S_l ) -\sum_{(j-1)r<k,l\le jr}\de(S_k,S_l
)\Big]+\E B_{n}-m\E B_r. \nonumber
\end{eqnarray}

We claim that
\begin{equation}
\E B_{n}-m\E B_r =O( n).\label{2.0w}
\end{equation} To see this, write
\begin{equation}
\E B_n-m\E B_r=\E B_n-\sum_{l=1}^m\E B(((l-1)r, lr]_<^2)\label{x2.3}
\end{equation} Notice that
\begin{eqnarray} && B_n-\sum_{l=1}^m B(((l-1)r, lr]_<^2)\label{x2.4}\\
&& =\sum_{l=1}^m B(((l-1)r, lr]\times (lr,mr])+B((mr, n]_<^2)\nn\\
&&\hspace{ 2in} +B((0, mr]\times (mr,n])\nn
\end{eqnarray} Since
\begin{equation} B(((l-1)r, lr]\times (lr,mr])\stackrel{d} {=}I_{ r,(
m-l)r}\label{x2.5}
\end{equation} by (\ref{2.0x}) we have
\begin{eqnarray} &&\E B(((l-1)r, lr]\times (lr,mr])\label{x2.1}\\ &&\le
C\Big\{(m-(l-1))r)\log(m-(l-1))r)\nn\\ && \hspace{
2in}-((m-l)r)\log((m-l)r) -r\log r\Big\}\nn
\end{eqnarray} Therefore
\begin{eqnarray} &&\sum_{l=1}^m\E  B(((l-1)r, lr]\times
(lr,mr])\label{x2.6}\\ &&\le
C\sum_{l=1}^m\Big\{(m-(l-1))r)\log(m-(l-1))r)\nn\\ && \hspace{
2in}-((m-l)r)\log((m-l)r) -r\log r\Big\}\nn\\ &&=C\Big\{mr\log
mr-mr\log r\Big\}=Cmr\log m.\nn
\end{eqnarray} Using (\ref{2.0y}) for $\E B((0, mr]\times (mr,n])=\E
I_{mr,i }$ and  (\ref{4.3}) for
$\E B((mr, n]_<^2)$ then completes the proof of (\ref{2.0w}).

Note that the summands in (\ref{2.30}) are independent. Therefore, for
some constant $C>0$ depending only on $\theta$ and $m$,
\begin{equation}
\E\exp\Big\{{\theta\over n}(\E B_{n}- B_{n})\Big\}\le C
\bigg(\E\exp\Big\{{\theta\over n}(\E B_r- B_r)\Big\}\bigg)^m
\label{2.31}
\end{equation} which proves (\ref{2.28}), since $\theta/n\leq
\theta/mr<\theta_o/r$ and
$r\rar\ff$ as $n\rar\ff$.

Then, by Chebychev's inequality, for any fixed $h>0$
\begin{equation}
\P\Big\{\E B_n -B_n
\ge \lambda nb_n\Big\}\leq e^{ -h\lambda b_n}\E\exp\Big\{{h\over
n}(\E B_n- B_n)\Big\}\label{2.32}
\end{equation} so that by (\ref{2.28})
\begin{equation}\lim_{n\to\infty}{1\over b_n}\log\P\Big\{\E B_n -B_n
\ge \lambda nb_n\Big\}\leq -h\lambda.
\label{2.33}
\end{equation} Since $h>0$ is arbitrary, this proves (\ref{3.3}).
\qed

\section{ Proof of Theorem \ref{theo-md}}

\medskip

By the G\"artner-Ellis theorem ( \cite[Theorem 2.3.6]{DZ}),   we need only
prove
\begin{equation}\hspace{ .3in}
\lim_{n\to\infty}{1\over b_n}\log\E\exp\Big\{
\theta\sqrt{b_n\over n}\vert B_n-\E B_n\vert^{1/2}\Big\} ={1\over
4}\kappa(2,2)^4\theta^2\det(\Gamma)^{-1/2}.
\label{3.1}
\end{equation}
      Indeed, by the G\"artner-Ellis theorem the above implies that
\begin{equation}
\lim_{n\to\infty}{1\over b_n}\log\P\Big\{\vert B_n-\E B_n\vert
\ge \lambda nb_n\Big\} =-\lambda\sqrt{\det(\Gamma)}\kappa (2,2)^{-4}.
\label{3.2}
\end{equation} Using (\ref{3.3}) we will then have Theorem
\ref{theo-md}. It thus  remains to prove (\ref{3.1}).

Let $f$ be a symmetric probability density function in the Schwartz space
$\mathcal{S}(\R^{2})$ of $C^{ \ff}$ rapidly decreasing functions.
      Let $\epsilon >0$ be a small number and write
\begin{equation}  f_\epsilon (x)=\epsilon^{-2}f(\epsilon^{-1}
x),\hskip.2in x\in\R^2
\label{3.6}
\end{equation} and
\begin{equation} l(n,x)=\sum_{k=1}^n \de( S_k,\,x),\hskip.2in
l(n,x,\epsilon)=\sum_{k=1}^n f_{\epsilon(b_n^{ -1}n)^{ 1/2}} (S_k-x).
\label{3.7}
\end{equation} By  \cite[Theorem 3.1]{C04},
\begin{eqnarray} &&\lim_{n\to\infty}{1\over b_n}\log\E\exp\bigg\{
{\theta\over\sqrt{2}}\sqrt{b_n\over
n}\Big(\sum_{x\in\Z^2}l^2(n,x,\epsilon)
\Big)^{1/2}\bigg\}
\label{3.8}\\ && =\sup_{g\in{\cal F}_2}\bigg\{{\theta\over\sqrt{2}}
\bigg(\int_{\R^2}\vert g^2\ast f_\epsilon
(x)\vert^2dx\bigg)^{1/2}-{1\over 2}\int_{\R^2}\langle \nabla g,\Gamma
\nabla g\rangle dx\bigg\}\nonumber
\end{eqnarray} where
\begin{equation} {\cal F}_2=\lc g\in W^{ 1,2}( \R^{ 2}) \,|\,\,\,\|g\|_{
2}=1\rc.\label{3.8a}
\end{equation}

As in the proof  of \cite[Theorem 1]{CR05}, (\ref{3.1}) will follow from
(\ref{3.8}) and the next Theorem.
\bt\label{theo-expapp}For any $\theta >0$,
\[\lim_{\ep\to 0}
\lim_{n\to\infty}{1\over b_n}\log\E\exp\Big\{
\theta\sqrt{b_n\over n}\vert B_n-\E B_n-{1 \over
2}\sum_{x\in\Z^2}l^2(n,x,\epsilon)\vert^{1/2}\Big\} =0.
\]
\et

\proof
\medskip Let $l>1$ be a large but fixed integer. Divide $[1 ,n]$ into
$l$ disjoint subintervals $D_1,\cdots, D_l$, each of length $[n/l]$ or
$[n/l]+1$.  Write
\begin{equation} D_i^*=\{(j,k)\in D^{ 2}_{ i};\hskip.1in j<k\}\hskip.2in
i=1,\cdots, l
\label{3.4}
\end{equation} With the notation of (\ref{2.9}) we have
\begin{equation}
      B_n=\sum_{i=1}^l  B(D_i^*)+\sum_{1\le j<k\le l} B(D_j\times D_k)
\label{3.5}
\end{equation}
Notice that
\begin{eqnarray} B(D_j\times D_k)&=&\sum_{ n_{ 1}\in D_j,n_{ 2}\in
D_k }\de(S_{ n_{ 1}},S_{ n_{ 2}} )
\nn\\ &=& \sum_{ n_{ 1}\in D_j,n_{ 2}\in D_k }\de((S_{ n_{ 1}}-S_{ b_{
j}})+S_{ b_{ j}},S_{ a_{ k}}+(S_{ n_{  2}}-S_{  a_{ k}}) )\nonumber
\\ &=& \sum_{ n_{ 1}\in D_j,n_{ 2}\in D_k }\de((S_{ n_{ 1}}-S_{  b_{
j}}),Z+(S_{ n_{ 2}}-S_{  a_{ k}}))\label{20.1}
\end{eqnarray}
      with $Z\buildrel d\over=S_{a_k}-S_{b_j}$, so that $Z, S_{ n_{ 1}}-S_{
b_{ j}},S_{ n_{ 2}}-S_{  a_{ k}}$ are  independent.  Then as in (\ref{2.1s})
\begin{equation}\qquad
\E  B(D_j\times D_k)=\E \sum_{ n_{ 1}\in D_j,n_{ 2}\in
D_k } p_{n_{ 1}+n_{ 2} }( Z)\leq
\sum_{ n_{ 1}\in D_j,n_{ 2}\in
D_k } p_{n_{ 1}+n_{ 2} }( 0)\label{20.1sup}
\end{equation}
since $\sup_{ x}p_{ j}( x)=p_{ j}( 0)$ for a symmetric random walk. Then
as in the proof of (\ref{2.2}) we have that
\begin{equation}\qquad
\E  B(D_j\times D_k)\leq cn/l.\label{20.1supl}
\end{equation}

Hence,
\begin{eqnarray} &&\hspace{ .2in}B_n-\E B_n
\label{3.5a}\\ &&=\sum_{i=1}^l
\big[ B(D_i^*)-\E  B(D_i^*)\big] +\sum_{1\le j<k\le l} B(D_j\times D_k)
-\E\sum_{1\le j<k\le l} B(D_j\times D_k) \nonumber\\ && =\sum_{i=1}^l
\big[ B(D_i^*)-\E  B(D_i^*)\big] +\sum_{1\le j<k\le l} B(D_j\times
D_k)+O(n)\nn
\end{eqnarray} where the last line follows from
(\ref{20.1supl}).

Write
\begin{equation}
\xi_i(n,x,\epsilon)=\sum_{k\in D_i} f_{\epsilon(b_n^{ -1}n)^{ 1/2}}
(S_k-x).
\label{3.9}
\end{equation} Then
\begin{equation}\hspace{ .4in}
\sum_{x\in\Z^2}l^2(n,x,\epsilon)
=\sum_{i=1}^l\sum_{x\in\Z^2}\xi_i^2(n,x,\epsilon) +2\sum_{1\le j\le
k\le l}\sum_{x\in\Z^2}\xi_j(n,x,\epsilon)\xi_k(n,x,\epsilon).
\label{3.10}
\end{equation} Therefore, by (\ref{3.5a})
\begin{eqnarray} &&\Big\vert (B_n-\E B_n)-{1\over 2}
\sum_{x\in\Z^2}l^2(n,x,\epsilon)\Big\vert
\label{3.11}\\ && \le\sum_{i=1}^l\big\vert  B(D_i^*)-\E
B(D_i^*)\big\vert + {1\over
2}\sum_{i=1}^l\sum_{x\in\Z^2}\xi_i^2(n,x,\epsilon)\nonumber\\ &&
+\sum_{1\le j<k\le l}\Big\vert  B(D_j\times D_k)
-\sum_{x\in\Z^2}\xi_j(n,x,\epsilon)\xi_k(n,x,\epsilon)\Big\vert +O(n).\nn
\end{eqnarray} The proof of Theorem \ref{theo-expapp} is completed in
the next two lemmas.

\bl\label{lem-ldbound} For any $\theta >0$,
\begin{eqnarray} &&
\limsup_{n\to\infty}{1\over b_n}\log\E\exp\Big\{
\theta\sqrt{b_n\over n}
\Big(\sum_{i=1}^l\sum_{x\in\Z^2}\xi_i^2(n,x,\epsilon)\Big)^{1/2}\Big\}
\label{3.12}\\&&\hspace{ .4in}\le l^{-1}{1\over
2}\kappa(2,2)^4\theta^2\det(\Gamma)^{-1/2}
\nn
\end{eqnarray} and
\begin{equation}\hspace{ .4in}
\limsup_{n\to\infty}{1\over b_n}\log\E\exp\Big\{
\theta\sqrt{b_n\over n}\Big(\sum_{i=1}^l\big\vert  B(D_i^*) -\E
      B(D_i^*)\big\vert\Big)^{1/2}\Big\}
\le l^{-1} H\theta^2,
\label{3.13}
\end{equation}
      where
\begin{equation}  H=\bigg(\sup\bigg\{\lambda >0;\hskip.1in
\sup_{ n}\E\exp\Big\{\lambda{1\over n}\vert B_n-\E
B_n\vert\Big\}<\infty
\bigg\}\bigg)^{-2}.
\label{3.14}
\end{equation}
\el

\proof  Replacing $\theta$ by $\theta/\sqrt l$, $n$ by $n/l$,
  and $b_n$ by $b^*_n=b_{ln}$
(notice that $b^*_{n/l}=b_n$)
\begin{eqnarray} &&
\limsup_{n\to\infty}{1\over b_n}
\log\E\exp\bigg\{\theta\sqrt{b_n\over n}
\bigg(\sum_{x\in\Z^2}\xi_i^2(n,x,\epsilon)\bigg)^{1/2}\bigg\}\label{3.14cle}\\ 
&&
=\limsup_{n\to\infty}{1\over b^*_{n/l}}
\log\E\exp\bigg\{{\theta \over \sqrt l}\sqrt{b^*_{n/l}\over n/l}
\bigg(\sum_{x\in\Z^2}\xi_i^2(n,x,\epsilon)\bigg)^{1/2}\bigg\} \nonumber
\end{eqnarray}

  Applying
Jensen's inequality on the right hand side of (\ref{3.8}),
\begin{eqnarray*}
\int_{\R^2} |g^2*f_\eps(x)|^2&=&\int_{\R^2}\Big[\int_{\R^2} 
g^2(x-y)f_\eps(y)\, dy\Big]^2\, dx\\
&\leq&\int \int g^4(x-y) f_\eps(y)\, dy\, dx=\int f_\eps(y)\Big[\int 
g^4(x-y)\, dx\Big]\, dy\\
&=&\Big[\int g^4(x)\, dx\Big]\int f_\eps(y)\, dy=\int_{\R^2} g^4(y)\, dy.
\end{eqnarray*}
Combining the last two displays with (\ref{3.8}) we have that
\begin{eqnarray} &&\limsup_{n\to\infty}{1\over b_n}
\log\E\exp\bigg\{\theta\sqrt{b_n\over n}
\bigg(\sum_{x\in\Z^2}\xi_i^2(n,x,\epsilon)\bigg)^{1/2}\bigg\}
\label{3.15}\\
&&\le \sup_{g\in {\cal F}_2}\bigg\{{\theta \over \sqrt l}
        \bigg(\int_{\R^2} \vert g(x)\vert^4dx\bigg)^{1/2}-{1\over
2}\int_{\R^2}
\langle \nabla g(x),\Gamma \nabla g(x)\rangle dx\bigg\}
\nonumber\\ && =l^{-1}\theta^2
\sup_{h\in {\cal F}_2}\bigg\{
        \bigg(\int_{\R^2} \vert h(x)\vert^4dx\bigg)^{1/2}-{1\over
2}\int_{\R^2}
\vert\nabla h(x)\vert^2 dx\bigg\}
      \nonumber\\ && ={1\over 2}l^{-1}\det (\Gamma)^{-1/2}\kappa(2,2)^4
\theta^2,\nn
\end{eqnarray}
       where the third line follows from the substitution
$g(x)=\sqrt{\vert\det (A)\vert}f(Ax)$ with a $2\times 2$ matrix $A$
satisfying
\begin{equation} A^{\tau}\Gamma A={1\over
2}\theta^2\sqrt{\det(\Gamma)}{\bf I}_2
\label{3.15m}
\end{equation}
      and the last line in \cite[Lemma A.2 ]{C04}; here ${\bf I}_2$ is the
$2\times 2$ identity matrix.
\medskip

Given $\delta >0$, there exist
$\ol{a}_1=(a_{1,1},\cdots, a_{1,l}),
\cdots, \ol{a}_m=(a_{m,1},\cdots, a_{m,l})$ in $\R^l$ such that
$\vert\ol{a}_1\vert =\cdots =\vert\ol{a}_m\vert =1$ and
\begin{equation}
\vert z\vert \le (1+\delta)\max\{\ol{a}_1\cdot z,\cdots,
\ol{a}_m\cdot z\}, \hskip.2in  z\in\R^l.
\label{3.16}
\end{equation}
       In particular, with
\begin{equation}
z=\bigg(\bigg(\sum_{x\in\Z^2}\xi_1^2(n,x,\epsilon)\bigg)^{1/2},\ldots,
\bigg(\sum_{x\in\Z^2}\xi_l^2(n,x,\epsilon)\bigg)^{1/2} \bigg)
\label{3.17}
\end{equation}
      we have
\begin{equation}\hspace{ .4in}
\bigg(\sum_{i=1}^l\sum_{x\in\Z^2}\xi_i^2(n,x,\epsilon)\bigg)^{1/2}
\le (1+\delta)\max_{1\le j\le m}\sum_{i=1}^l a_{j,i}
\bigg(\sum_{x\in\Z^2}\xi_i^2(n,x,\epsilon)\bigg)^{1/2}.
\label{3.18}
\end{equation}
       Hence
\begin{eqnarray} &&
\E\exp\bigg\{\theta\sqrt{b_n\over n}
\bigg(\sum_{i=1}^l\sum_{x\in\Z^2}\xi_i^2(n,x,\epsilon)\bigg)^{1/2}\bigg\}
\label{3.19}\\ &&
\le\sum_{j=1}^m\E\exp\bigg\{\theta\sqrt{b_n\over n}
(1+\delta)\sum_{i=1}^l a_{j,i}
\bigg(\sum_{x\in\Z^2}\xi_i^2(n,x,\epsilon)\bigg)^{1/2}\bigg\}\nonumber\\
&&
      =
\sum_{j=1}^m\prod_{i=1}^l
\E\exp\bigg\{\theta\sqrt{b_n\over n} (1+\delta)a_{j,i}
\bigg(\sum_{x\in\Z^2}\xi_i^2(n,x,\epsilon)\bigg)^{1/2}\bigg\},\nonumber
\end{eqnarray}
      where the last line follows from independence of
$\| \xi_i(n,x,\epsilon)  \|_{ L^{ 2}(\Z^2 )}$, $i=1,
\ldots,l$.  Therefore
\begin{eqnarray} &&
\limsup_{n\to\infty}{1\over b_n}
\log\E\exp\bigg\{\theta\sqrt{b_n\over n}
\bigg(\sum_{k=1}^l\sum_{x\in\Z^2}\xi_k^2(n,x,\epsilon)\bigg)^{1/2}\bigg\}
\label{3.20}\\ && \le\max_{1\le j\le m}{1\over 2}l^{-1}\kappa(2,2)^4
(1+\delta)^2\theta^2\Big(\sum_{i=1}^la_{j,i}^2\Big)\nonumber
\\ && ={1\over 2}l^{-1}\det(\Gamma)^{-1/2}
\kappa(2,2)^4(1+\delta)^2\theta^2.\nonumber
\end{eqnarray}
       Letting $\delta\to 0^+$ proves (\ref{3.12}).

\medskip

By the inequality $ab\leq a^{ 2}+b^{ 2}$ we have that
\begin{eqnarray} &&
\E\exp\Big\{\theta\sqrt{b_n\over n}\vert B_n-\E B_n\vert^{1/2}\Big\}
\label{3.21}\\&&
\le \exp\big\{c^2\theta^2b_n\big\}\E\exp\Big\{c^{-2}{1\over n}\vert
B_n-\E B_n\vert\Big\},
\nn
\end{eqnarray} and taking $c^{ -2}\uparrow H^{ -2}$ we see that for any
$\theta >0$,
\begin{equation}
\limsup_{n\to\infty}{1\over b_n}\log
\E\exp\Big\{\theta\sqrt{b_n\over n}\vert B_n-\E B_n\vert^{1/2}\Big\}
\le H^2\theta^2.
\label{3.22}
\end{equation}
      Notice that for any $1\le i\le l$,
\begin{equation}
      B(D_i^*)-\E B(D_i^*)\buildrel d\over = B_{\#(D_i)}-\E B_{\#(D_i)}.
\label{3.23}
\end{equation}
  We have
\[
\E\exp\Big\{\theta\sqrt{b_n\over n}\vert B(D_i^*)-\E
B(D_i^*)\vert^{1/2}\Big\}=
\E \exp\Bigl\{ \frac{\theta}{\sqrt l} 
\sqrt{\frac{b_n}{n/l}}|B(D_i^*-\E B(D_i^*)|\Bigr\}.
\]
Replacing $\theta$ by $\theta/\sqrt l$, $n$ by $n/l$, and $b_n$ by 
$b^*_n=b_{ln}$
(notice that $b^*_{n/l}=b_n$)
gives
\begin{equation}
\limsup_{n\to\infty}{1\over b_n}\log
\E\exp\Big\{\theta\sqrt{b_n\over n}\vert B(D_i^*)-\E
B(D_i^*)\vert^{1/2}\Big\}
\le l^{-1}H^2\theta^2.
\label{3.24}
\end{equation} Thus (\ref{3.13}) follows by the same argument we used to
prove  (\ref{3.12}).\qed

\bl\label{lem-ldapprox} For any $\theta >0$ and any $1\le j <k\le l$,
\begin{equation}
\limsup_{\epsilon\to 0^+}\limsup_{n\to\infty}{1\over
b_n}\log\E\exp\Big\{
\theta\sqrt{b_n\over n}\big\vert  B(D_j\times
D_k)-\sum_{x\in\Z^2}\xi_j(n,x,\epsilon)\xi_k(n,x,\epsilon)
\big\vert^{1/2}\Big\}=0.
\label{3.25}
\end{equation}
\el

\proof Define $a_j,b_j$ so that $D_j=(a_j,b_j]$ $(1\le j\le l)$. We now fix
$1\le j <k\le l$ and estimate
\begin{equation} B(D_j\times
D_k)-\sum_{x\in\Z^2}\xi_j(n,x,\epsilon)\xi_k(n,x,\epsilon).
\label{3.26}
\end{equation}
      Without loss of generality we may assume that
$v=:[n/l]=\#(D_j)=\#(D_k)$.  For $y\in Z^{ 2}$ set
\begin{equation} I_{n}( y)=\sum_{ n_{ 1},n_{ 2}=1}^{ n}
\de( S_{ n_{ 1}},S'_{ n_{ 2}}+y).\label{18.1}
\end{equation} Note that
$I_{ n}=I_{n}(0)$. By (\ref{20.1}) we have that
\begin{equation} B(D_j\times D_k)\stackrel{d}{=}I_{v}( Z)\label{20.2}
\end{equation} with $Z$ independent of $S,S'$.

Similarly, we have
\begin{eqnarray}&&
\sum_{x\in\Z^2}\xi_j(n,x,\epsilon)\xi_k(n,x,\epsilon)\nn\\
&&=\sum_{x\in\Z^2}\sum_{ n_{ 1}\in D_k,n_{ 2}\in D_k
}f_{\epsilon(b_n^{ -1}n)^{ 1/2}} (S_{ n_{ 1}}-x)f_{\epsilon(b_n^{ -1}n)^{
1/2}} (S_{ n_{ 2}}-x)
\nn\\ &&=\sum_{x\in\Z^2}\sum_{ n_{ 1}\in D_k,n_{ 2}\in D_k
}f_{\epsilon(b_n^{ -1}n)^{ 1/2}} (x)f_{\epsilon(b_n^{ -1}n)^{ 1/2}} (S_{
n_{ 2}}-S_{ n_{ 1}}-x)
\nn\\&& =\sum_{ n_{ 1}\in D_k,n_{ 2}\in D_k }f_{\epsilon(b_n^{ -1}n)^{
1/2}}\circledast f_{\epsilon(b_n^{ -1}n)^{ 1/2}} (S_{ n_{ 2}}-S_{ n_{ 1}})
\nn\\&& =\sum_{ n_{ 1}\in D_j,n_{ 2}\in D_k }f_{\epsilon(b_n^{ -1}n)^{
1/2}}\circledast f_{\epsilon(b_n^{ -1}n)^{ 1/2}}((S_{ n_{ 2}}-S_{  a_{
k}})-(S_{ n_{ 1}}-S_{  b_{ j}})+Z)\label{20.1a}
\end{eqnarray} where
\begin{equation} f\circledast f( y)=
\sum_{x\in\Z^2}f( x)f( y-x)\label{20.0}
\end{equation} denotes convolution  in $ L^{ 1}(\Z^{ 2})$. It is clear that if
$f\in \mathcal{S}(\R^{2})$ so is $f\circledast f$. For $y\in Z^{ 2}$, define
the link
\be L_{ n,\ep}( y)=\sum_{ n_{ 1},n_{ 2}=1}^{ n} f_{\epsilon}\circledast
f_{\epsilon}(S'_{ n_{ 2}}- S_{ n_{ 1}}+y).\label{18.3}
\ee By (\ref{20.1a}) we have that
\begin{equation}
\sum_{x\in\Z^2}\xi_j(n,x,\epsilon)\xi_k(n,x,\epsilon)\stackrel{d}{=}
L_{v,(b_n^{-1}n)^{ 1/2}\ep}( Z)\label{20.2a}
\end{equation} with $Z$ independent of $S,S'$.

\bl\label{lem-conv}Let $f\in \mathcal{S}(\R^{2})$ with Fourier transform
$\wh{f}$ supported on $(-\pi,\pi)^{ 2}$. Then for any $r\geq 1$
\begin{equation}
\int e^{ -i\la y}(f_{r}\circledast f_{r})( y)\,dy=\wh{f} ^{ 2}( r\la),\hspace{
.2in}\forall \la\in \R^{ 2}.\label{20.0dd}
\end{equation}
\el

\proof We have
\begin{eqnarray}
\int e^{ -i\la y}(f\circledast f)( y)\,dy&=&
\sum_{x\in\Z^2}f( x)\,\int e^{ -i\la y}f( y-x)\,dy\label{20.0a}\\ &=&
\wh{f}( \la)\sum_{x\in\Z^2}f( x)e^{ -i\la x} \nonumber\\ &=&
\wh{f}( \la)\sum_{x\in\Z^2}\(\int e^{ ip x}\wh{f}( p)\,dp\)e^{ -i\la x}.
\nonumber
\end{eqnarray} For $x\in\Z^2$
\begin{eqnarray} &&
\int e^{ ip x}\wh{f}( p)\,dp=\sum_{u\in\Z^2}\int_{[-\pi,\pi]^{ 2} } e^{ ip
x}\wh{f}( p+2\pi u)\,dp\label{20.0b}\\ && \nonumber
\end{eqnarray} and using Fourier inversion
\begin{eqnarray} &&
\sum_{x\in\Z^2}\(\int e^{ ip x}\wh{f}( p)\,dp\)e^{ -i\la x}\nn\\ &&
=\sum_{u\in\Z^2}\sum_{x\in\Z^2}\(
\int_{[-\pi,\pi]^{ 2} } e^{ ip x}\wh{f}( p+2\pi u)\,dp\)e^{ -i\la
x}\label{20.0bb}\\ &&=\sum_{u\in\Z^2}\wh{f}(\la+2\pi u).\nn
\end{eqnarray} Thus from (\ref{20.0a}) we find that
\begin{equation}
\int e^{ -i\la y}f\circledast f( y)\,dy=\wh{f}(
\la)\sum_{u\in\Z^2}\wh{f}(\la+2\pi u).\label{20.0c}
\end{equation} Since $\wh{f}_{ r}( \la)=\wh{f}( r\la)$ we see that for any
$r>0$
\begin{equation}
\int e^{ -i\la y}(f_{r}\circledast f_{r})( y)\,dy=\wh{f}(
r\la)\sum_{u\in\Z^2}\wh{f}(r\la+2\pi ru).\label{20.0d}
\end{equation} Then if  $r\geq1$, using the fact that $\wh{f}(
\la)$ is supported in
$(-\pi,\pi)^{ 2}$, we obtain (\ref{20.0dd}).\qed

Taking $f\in {\cal S}( R^2)$ with $\wh{f}( \la)$ supported in
$(-\pi,\pi)^{ 2}$,
      Lemma \ref{lem-ldapprox} will follow from Theorem
\ref{theo-exprw} of the next section.
\qed

\section{Intersections of Random Walks}\label{sec-rwexpapp}

      Let $S_1(n), S_2(n)$  be independent copies of the symmetric random
walk
$S(n)$ in $Z^{ 2}$ with a finite second moment.

Let
$f$ be a positive symmetric function in the Schwartz space
${\cal S}( R^2)$ with $\int f\,dx=1$ and $\wh f$ supported in
$(-\pi,\pi)^{ 2}$. Given $\epsilon >0$, and with the notation of the last
section, let us define the
     link
\be I_{ n,\ep}( y)=\sum_{ n_{ 1},n_{ 2}=1}^{ n} f_{(b_n^{ -1}n)^{
1/2}\epsilon}\circledast f_{(b_n^{ -1}n)^{ 1/2}\epsilon} (S_2(n_2)-
S_1(n_1)+y))\label{18.3a}
\ee with $I_{ n,\ep}=I_{n,\ep}(0)$.
\medskip

\bt\label{theo-exprw} For any $\la>0$
\begin{eqnarray} &&
\limsup_{ \ep\rar 0}\,\limsup_{ n\rar
\ff}\,\sup_{y}\label{18.20}\\ &&{ 1\over b_n}\log E\(\exp\lc
\la \Bigg|  {I_{ n}( y)-I_{n,\ep}( y)
\over b_n^{ -1}n}\Bigg|^{1/2}\rc\)=0.\nn
\end{eqnarray}
\et

\noindent{\bf Proof of Theorem \ref{theo-exprw}.}  We have
\begin{eqnarray} &&  {1 \over b_n^{ -1}n} I_{ n}(y)\label{18.7}\\ && ={1
\over b_n^{ -1}n}\sum_{ n_{ 1},n_{ 2}=1}^{ n} \de(
S_1(n_1),S_2(n_2)+y)\nn\\  &&  ={1
\over b_n^{ -1}n( 2\pi )^{2}}\sum_{ n_{ 1},n_{ 2}=1}^{ n}\bigg[
      \int_{ [-\pi,\pi]^{2}}e^{ i p\cdot ( S_2(n_2)+y-S_1(n_1))}
\, dp\bigg]\nonumber
\end{eqnarray}  where from now on we work modulo $\pm
\pi$. Then by scaling we have
\begin{eqnarray} && {1 \over b_n^{ -1}n} I_{ n}(y)\label{18.8}\\ &&
      ={1 \over (b_n^{ -1}n)^{ 2}( 2\pi )^{2}}\sum_{ n_{ 1},n_{ 2}=1}^{
n}\bigg[
      \int_{(b_n^{ -1} n)^{ 1/2} [-\pi,\pi]^{ 2}} e^{  i p\cdot (
S_2(n_2)+y-S_1(n_1))/(b_n^{ -1} n)^{ 1/2}}
\,  dp\bigg]\nonumber
\end{eqnarray}

As in (\ref{18.7})-(\ref{18.8}), using Lemma \ref{lem-conv}, the fact that
    $\ep (b_n^{ -1} n)^{ 1/2}\geq 1$ for $\ep>0$ fixed and large enough
$n$, and abbreviating $\wh{h}=\wh{f} ^{ 2}$
\begin{eqnarray}  &&  {1 \over b_n^{ -1}n} I_{ n,\ep}( y)\label{18.8a}\\
      &&  ={1 \over b_n^{ -1}n( 2\pi )^{2}}\sum_{ n_{ 1},n_{ 2}=1}^{ n}\bigg[
\int_{ R^{ 2}}e^{ i p\cdot ( S_2(n_2)+y-S_1(n_1))}
\, \wh{h}( \ep (b_n^{ -1} n)^{ 1/2}p)\, dp\bigg]\nn\\
      &&  = {1 \over (b_n^{ -1}n)^{ 2}( 2\pi )^{2}}\sum_{ n_{ 1},n_{ 2}=1}^{
n}\bigg[
      \int_{R^{ 2}} e^{ i p\cdot
      (S_2(n_2)+y-S_1(n_1))/(b_n^{ -1} n)^{ 1/2}}
\, \wh{h}( \ep p)\,  dp\bigg].\nonumber
\end{eqnarray} Using our assumption that  $\wh h$ supported in
$[-\pi,\pi]^{ 2}$, and that $\ep^{ -1}\leq (b_n^{ -1} n)^{ 1/2}$ for
$\ep>0$ fixed and large enough $n$, we have that
\begin{eqnarray}  &&  {1 \over b_n^{ -1}n} I_{ n,\ep}( y)\label{18.8ay}\\
      &&  = {1 \over (b_n^{ -1}n)^{ 2}( 2\pi )^{2}}\sum_{ n_{ 1},n_{ 2}=1}^{
n}\nn\\
      &&\hspace{ .3in}\bigg[
      \int_{\ep^{ -1} [-\pi,\pi]^{2}} e^{ i p\cdot
      (S_2(n_2)+y-S_1(n_1))/(b_n^{ -1} n)^{ 1/2}}
\, \wh{h}( \ep p)\,  dp\bigg]\nonumber\\
      &&  = {1 \over (b_n^{ -1}n)^{ 2}( 2\pi )^{2}}\sum_{ n_{ 1},n_{ 2}=1}^{
n}\nn\\
      &&\hspace{ .3in}\bigg[
      \int_{(b_n^{ -1} n)^{ 1/2} [-\pi,\pi]^{2}} e^{ i p\cdot
      (S_2(n_2)+y-S_1(n_1))/(b_n^{ -1} n)^{ 1/2}}
\, \wh{h}( \ep p)\,  dp\bigg].\nonumber
\end{eqnarray}

To prove (\ref{18.20}) it suffices to  show that for each
$\la>0$ we have
\begin{eqnarray} && \sup_{y}E\(\exp \lc \la
\Bigg|{I_{ n}( y)-I_{ n,\ep}( y) \over b_{ n}^{ -1}n}\Bigg|^{ 1/2}\rc\)
\label{18.16c2}\\ &&
\leq Cb_n( 1-C\la \ep^{m/4 })^{ -1}( 1+C\la \ep^{1/4}b_n^{
1/2})e^{C\la^{ 2}\ep^{1/2}b_n}.\nn
\end{eqnarray}
      for some $C<\ff$ and all $\ep>0$ sufficiently small.

      We begin by expanding
\begin{eqnarray} && E\(\exp \lc \la
\Bigg|{I_{ n}( y)-I_{ n,\ep}( y) \over b_{ n}^{ -1}n}\Bigg|^{ 1/2}\rc\)
\label{18.16cm}\\ &&=\sum_{ m=0}^{ \ff}{\la^{ m} \over m!} E\(\Bigg|  {1
\over b_n^{ -1}n}(  I_{ n}( y)-I_{n,\ep}( y))\Bigg|^{ m/2}\)\nn\\ &&\leq
\sum_{ m=0}^{ \ff}{\la^{ m} \over m!}
\(E\(\lc  {1
\over b_n^{ -1}n}(  I_{ n}( y)-I_{n,\ep}( y))\rc^{ 2m}\)\)^{ 1/4}\nn
\end{eqnarray}

       By (\ref{18.8}), (\ref{18.8ay}) and the symmetry of $S_{1}$ we have
\begin{eqnarray} && E\(\lc  {1
\over b_n^{ -1}n}(  I_{ n}( y)-I_{n,\ep}( y))\rc^{ m}
\)\label{18.13ch}\\ &&
      ={1 \over (b_n^{ -1}n)^{ 2m}( 2\pi )^{2m}}
\sum_{\stackrel{n_{ 1,j},n_{ 2,j}=1}{j=1,\ldots,m}}^{n}
\int_{(b_n^{ -1} n)^{ 1/2} [-\pi,\pi]^{ 2m}}\nn\\
      &&\hspace{.5in} E\(e^{ i\sum_{ j=1}^{ m}p_{ j}\cdot
      (S_2(n_{2,j})+y+S_1(n_{1,j})) /(b_n^{ -1} n)^{ 1/2}} \)
\prod_{ j=1}^{m} (1-\wh{h}(\ep p_{ j} )) \,dp_{ j}.\nn
\end{eqnarray} Then
\begin{eqnarray} &&\qquad \Bigg |E\(\lc  {1
\over b_n^{ -1}n}(  I_{ n}( y)-I_{n,\ep}( y))\rc^{ m}
\)\Bigg |\label{18.13chswp}\\ &&
\leq {1 \over (b_n^{ -1}n)^{ 2m}( 2\pi )^{2m}}
      \sum_{\stackrel{n_{ 1,j} =1}{j=1,\ldots,m}}^{n}\sum_{\stackrel{n_{
2,j}=1}{j=1,\ldots,m}}^{n}
\int_{(b_n^{ -1} n)^{ 1/2} [-\pi,\pi]^{ 2m}}\nn\\
      &&\hspace{.5in}\Bigg |E\(e^{ i\sum_{ j=1}^{ m}p_{ j}\cdot S_1(n_{1,j})
/(b_n^{ -1} n)^{ 1/2}} \)\Bigg |
\nn\\
      &&\hspace{.5in}
      \Bigg |E\(e^{ i\sum_{ j=1}^{ m}p_{ j}\cdot
      S_2(n_{2,j}) /(b_n^{ -1} n)^{ 1/2}} \)\Bigg |
\prod_{ j=1}^{m} |1-\wh{h}(\ep p_{ j} )| \,dp_{ j}.\nn
\end{eqnarray}

By the Cauchy-Schwarz inequality
\begin{eqnarray} &&
\int_{(b_n^{ -1} n)^{ 1/2} [-\pi,\pi]^{ 2m}}
\Bigg |E\(e^{ i\sum_{ j=1}^{ m}p_{ j}\cdot S_1(n_{1,j}) /(b_n^{ -1} n)^{
1/2}} \)\Bigg |
\label{18.13chswbass}\\
      &&\hspace{.5in}
      \Bigg |E\(e^{ i\sum_{ j=1}^{ m}p_{ j}\cdot
      S_2(n_{2,j}) /(b_n^{ -1} n)^{ 1/2}} \)\Bigg |
\prod_{ j=1}^{m} |1-\wh{h}(\ep p_{ j} )| \,dp_{ j}\nn\\ &&
\leq \prod_{ i=1}^{ 2}\lc\int_{(b_n^{ -1} n)^{ 1/2} [-\pi,\pi]^{
2m}}\right.\nn\\
      &&\hspace{.4in} \left.\Bigg |E\(e^{ i\sum_{ j=1}^{ m}p_{ j}\cdot
S(n_{i,j}) /(b_n^{ -1} n)^{ 1/2}} \)\Bigg |^{ 2}
\prod_{ j=1}^{m} |1-\wh{h}(\ep p_{ j} )|  \,dp_{ j}\rc^{ 1/2}.\\ &&
\nonumber
\end{eqnarray} Thus
\begin{eqnarray} &&\qquad \Bigg |E\(\lc  {1
\over b_n^{ -1}n}(  I_{ n}( y)-I_{n,\ep}( y))\rc^{ m}
\)\Bigg |^{ 1/2}\label{18.13chsw}\\ &&
\leq \sum_{\stackrel{n_{ j} =1}{j=1,\ldots,m}}^{n}{1 \over (b_n^{ -1}n)^{
m}( 2\pi )^{m}}
\lc \int_{(b_n^{ -1} n)^{ 1/2} [-\pi,\pi]^{ 2m}}\right.\nn\\
      &&\hspace{.4in} \left.\Bigg |E\(e^{ i\sum_{ j=1}^{ m}p_{ j}\cdot
S(n_{j}) /(b_n^{ -1} n)^{ 1/2}} \)\Bigg |^{ 2}
\prod_{ j=1}^{m} |1-\wh{h}(\ep p_{ j} )|  \,dp_{ j}\rc^{ 1/2}\nn
\end{eqnarray}

For any permutation
$\pi$
      of $\{1.\ldots,m\}$ let
\begin{equation}\qquad D_{m}(\pi)=\{  (n_{1},\ldots,n_{ m})|\, 1\leq
n_{\pi( 1)}\leq \cdots \leq n_{\pi( m)}\leq n\}.\label{18.13j}
\end{equation}
       Using the (non-disjoint) decomposition
      \[\{ 1,\ldots,n\}^{ m}=\bigcup_{ \pi} D_{ m}(\pi)\] we have from
(\ref{18.13chsw}) that
\begin{eqnarray} && \Bigg | E\(\lc  {1
\over b_n^{ -1}n}(  I_{ n}( y)-I_{n,\ep}( y))\rc^{ m}
\)\Bigg |^{ 1/2}\label{18.13}\\ &&
\leq \sum_{ \pi} \sum_{ D_{ m}(\pi)} {1 \over (b_n^{ -1}n)^{ m}( 2\pi
)^{m}}
\lc\int_{(b_n^{ -1} n)^{ 1/2} [-\pi,\pi]^{2m}}\right.\nn\\
      &&\hspace{.5in}\left.\Bigg | E\(e^{ i\sum_{ j=1}^{ m}p_{ j}\cdot
      S(n_{j}) /(b_n^{ -1} n)^{ 1/2}} \)\Bigg |^{ 2}
\prod_{ j=1}^{m} |1-\wh{h}(\ep p_{ j} )| \,dp_{ j}\rc^{ 1/2}.\nn
\end{eqnarray} where the first sum is over all permutations
$\pi$
      of $\{1.\ldots,m\}$.

Set
\begin{equation}
\phi(u)=E\(e^{ i u\cdot S(1 )}\).\label{18.10}
\end{equation} It follows from our assumptions that
$\phi(u)\in C^{ 2}$, $ {\partial \over \partial u_{ i}}\phi(0)=0$ and
${\partial^{ 2} \over \partial u_{ i}\partial u_{ j}}\phi(0)=-E\( S_{ (i)}(1
)S_{(j)}(1 )\)$ where
$S(1 )=(  S_{ (1)}(1 ),S_{(2)}(1 ))$ so that for some
$\de>0$
\begin{equation}
\phi(u)=1-E\((u\cdot S(1 ))^{ 2}\)/2+o(|u|^{ 2} ),\hspace{ .2in}|u|\leq
\de.\label{18.10Tay}
\end{equation} Then for some $c_{ 1}>0$
\begin{equation}
\phi(u)\leq e^{ -c_{ 1}|u|^{ 2}},\hspace{ .2in}|u|\leq
\de.\label{18.10Tay2a}
\end{equation} Strong aperiodicity implies that $|\phi(u)|< 1$ for
$u\neq 0$ and
$u\in [-\pi,\pi]^{ 2}$. In particular, we can find $b<1$ such that
$|\phi(u)|\leq b$ for
$\de\leq |u|$ and $ u\in [-\pi,\pi]^{ 2}$. But clearly we can choose
$c_{ 2}>0$ so that
$b\leq e^{ -c_{ 2}|u|^{ 2}}$ for $u\in [-\pi,\pi]^{ 2}$. Setting
$c=\min(c_{ 1}, c_{ 2})>0$ we then have
\begin{equation}
\phi(u)\leq e^{ -c|u|^{ 2}},\hspace{ .2in}u\in [-\pi,\pi]^{
2}.\label{18.10Tay2}
\end{equation}

On
$D_{ m}(\pi)$ we can write
\begin{equation}\qquad
\sum_{ j=1}^{ m}p_{ j}\cdot S(n_{ j} ) =\sum_{ j=1}^{ m}(\sum_{ i=j}^{ m}
p_{\pi( i)})(S(n_{ \pi( j)} ) -S(n_{ \pi( j-1)} )).\label{18.14}
\end{equation}
      Hence on $D_{ m}(\pi)$
\begin{equation}
\qquad E\(e^{ i\sum_{ j=1}^{ m}p_{ j}\cdot
      S(n_{j}) /(b_n^{ -1} n)^{ 1/2}} \)= \prod_{ j=1}^{ m} \phi((\sum_{
i=j}^{ m} p_{\pi( i)}) /(b_n^{ -1} n)^{1/2})^ {(n_{\pi( j)}-n_{\pi(
j-1)})}.\label{18.15}
\end{equation} Now it is clear that
\begin{eqnarray} &&\sum_{ D_{ m}(\pi)}
\lc\int_{(b_n^{ -1} n)^{ 1/2} [-\pi,\pi]^{2m}}\right.\label{18.14x}\\
&&\hspace{ .2in}\left.\Bigg | \prod_{ j=1}^{ m} \phi((\sum_{ i=j}^{ m}
p_{\pi( i)}) /(b_n^{ -1} n)^{1/2})^ {(n_{\pi( j)}-n_{\pi( j-1)})}\Bigg |^{ 2}
\prod_{ j=1}^{m} |1-\wh{h}(\ep p_{ j} )| \,dp_{ j}\rc^{ 1/2}\nn
\\ &&=\sum_{ 1\leq n_{\pi( 1)}\leq \cdots \leq n_{\pi( m)}\leq n}
\lc\int_{(b_n^{ -1} n)^{ 1/2} [-\pi,\pi]^{2m}}\right.\nn\\ &&\hspace{
.2in}\left.\Bigg | \prod_{ j=1}^{ m} \phi((\sum_{ i=j}^{ m} p_{\pi( i)})
/(b_n^{ -1} n)^{1/2})^ {(n_{\pi( j)}-n_{\pi( j-1)})}\Bigg |^{ 2}
\prod_{ j=1}^{m} |1-\wh{h}(\ep p_{ j} )| \,dp_{ j}\rc^{ 1/2}\nn
\end{eqnarray} is independent of the permutation $\pi$. Hence writing
\begin{equation} u_{ j}=\sum_{ i=j}^{ m} p_{u}\label{18.14a}
\end{equation} we have from (\ref{18.13}) that
\begin{eqnarray} && \Bigg | E\(\lc  {1
\over b_n^{ -1}n}(  I_{ n}( y)-I_{n,\ep}( y))\rc^{ m}
\)\Bigg |^{ 1/2}\label{18.13rxy}\\ &&
\leq m! \sum_{ 1\leq n_{1}\leq \cdots \leq n_{m}\leq n} {1 \over (b_n^{
-1}n)^{ m}( 2\pi )^{m}}
\lc\int_{(b_n^{ -1} n)^{ 1/2} [-\pi,\pi]^{2m}}\right.\nn\\
      &&\hspace{.5in}\left.\Bigg | \prod_{ j=1}^{ m} \phi(u_{ j} /(b_n^{ -1}
n)^{1/2})^ {(n_{j}-n_{j-1})}\Bigg |^{ 2}
\prod_{ j=1}^{m} |1-\wh{h}(\ep p_{ j} )| \,dp_{ j}\rc^{ 1/2}.\nn
\end{eqnarray}

For each $A\subseteq \{ 2,3,\ldots,m\}$ we use $D_{ m}(  A)$ to denote
the subset of
$\{ 1\leq n_{1}\leq \cdots \leq n_{m}\leq n\}$ for which $n_{j}=n_{ j-1}$
if and only if $j\in A$. Then we have
\begin{eqnarray} && \Bigg | E\(\lc  {1
\over b_n^{ -1}n}(  I_{ n}( y)-I_{n,\ep}( y))\rc^{ m}
\)\Bigg |^{ 1/2}\label{18.13rx}\\ &&
\leq m! \sum_{A\subseteq \{ 2,3,\ldots,m\} }\sum_{D_{ m}(  A)} {1
\over (b_n^{ -1}n)^{ m}( 2\pi )^{m}}
\lc\int_{(b_n^{ -1} n)^{ 1/2} [-\pi,\pi]^{2m}}\right.\nn\\
      &&\hspace{.5in}\left.\Bigg | \prod_{ j=1}^{ m} \phi(u_{ j} /(b_n^{ -1}
n)^{1/2})^ {(n_{j}-n_{j-1})}\Bigg |^{ 2}
\prod_{ j=1}^{m} |1-\wh{h}(\ep p_{ j} )| \,dp_{ j}\rc^{ 1/2}.\nn
\end{eqnarray}

For any $u\in R^{ d}$ let $\wt{u}$ denote the representative of
$u\mbox{ mod }\,(b_n^{ -1} n)^{ 1/2}2\pi Z^{ 2}$ of smallest absolute
value. We note that
\begin{equation} |\wt{-u}|=|\wt{u}|,\hspace{ .2in}\mbox{ and }\hspace{
.2in} |\wt{u+v}|=|\wt{u}+\wt{v}|\leq |\wt{u}|+|\wt{v}|.\label{18.15obs}
\end{equation} Using the periodicity of $\phi$ we see that
(\ref{18.10Tay2}) implies that for all $u$
\begin{equation} |\phi(u/(b_n^{ -1} n)^{ 1/2})|
\leq e^{-c|\wt{u}|^{2}/(b_n^{ -1}n) }.\label{18.16av}
\end{equation} Then we have that on $\{ 1\leq n_{1}\leq \cdots \leq
n_{m}\leq n\}$
\begin{equation}\qquad\Bigg | \prod_{ j=1}^{ m} \phi(u_{ j} /(b_n^{ -1}
n)^{1/2})^ {(n_{j}-n_{j-1})}\Bigg |^{ 2}\leq \prod_{ j=1}^{ m}
e^{-c|\wt{u}_{ j}|^{ 2}(n_{j}-n_{j-1})/(b_n^{-1}n) } \label{x18.15a}
\end{equation}

Using  $| 1-\wh{h}(\ep p_{ j} )  |\leq c\ep^{1/2}p_{ j}^{ 1/2}$ we bound
the integral in (\ref{18.13rx}) by
\begin{equation} c^{ m}\ep^{ m/2}\int_{(b_n^{ -1} n)^{ 1/2}
[-\pi,\pi]^{2m}}
      \prod_{ j=1}^{ m} e^{-c|\wt{u}_{ j}|^{ 2}(n_{j}-n_{j-1})/(b_n^{-1}n) }
      |p_{ j}|^{1/2 }
\,dp_{ j}.\label{x6.6j}
\ee

      Using (\ref{18.14a}) and (\ref{18.15obs}) we have that
\begin{equation}
      \prod_{ j=1}^{ m} |p_{ j}|^{1/2 }\leq   \prod_{ j=1}^{
m}(|\wt{u}_{j}|^{1/2 }+|\wt{u}_{j+1}|^{1/2 })\label{x6.6jm}
\end{equation} and when we expand the right hand side as a sum of
monomials we can be sure that no factor
$|\wt{u}_{k}|^{1/2 }$ appears more than twice. Thus we see that we can
bound (\ref{x6.6j}) by
\begin{equation} \hspace{ .4in}C^{ m}\ep^{   m/2}\max_{h(j)
}\int_{(b_n^{ -1} n)^{ 1/2} [-\pi,\pi]^{2m}} \prod_{ j=1}^{ m}
e^{-c|\wt{u}_{ j}|^{ 2}(n_{j}-n_{j-1})/(b_n^{-1}n) }
      |\wt{u}_{ j}|^{ h(j)/2}
\,dp_{ j}\label{x6.6k}
\end{equation} where the $\max$ runs over the the set of functions
$h( j)$ taking values
$0,1\mbox{ or }2$ and such that $\sum_{j}h(j)=m$. Changing variables,
we thus need to bound
\begin{equation}\quad
\int_{ \La_{ n}}
\prod_{ j=1}^{ m} e^{-c|\wt{u}_{ j}|^{ 2}(n_{j}-n_{j-1})/(b_n^{-1}n) }
      |\wt{u}_{ j}|^{h(j)/2}
\,du_{ j}\label{x6.6l}
\end{equation} where, see (\ref{18.14a}),
\begin{equation}\qquad
\La_{ n}=\{ (u_{ 1},\ldots,u_{ m} )\,|\, u_{ j}-u_{ j+1}\in (b_n^{  -1} n)^{
1/2} [-\pi,\pi]^{2},\,\forall j\}.\label{x6.6m}
\end{equation}

Let $C_{ n}$ denote the rectangle $(b_n^{ -1} n)^{ 1/2} [-\pi,\pi]^{2}$
and let us call any rectangle of the form
$2\pi k+C_{ n}$, where $k\in Z^{ 2}$, an elementary rectangle. Note that
any rectangle of the form $v+C_{ n}$, where $v\in R^{ 2}$, can be covered
by $4$ elementary rectangles. Hence for any $v\in R^{2}$ and
$1\leq s\leq n$
\bea &&
\int_{v+C_{ n}}e^{-c{ s \over b_n^{ -1}n}|\wt{u}|^{2} }|\wt{u}|^{
a}\,du\label{18.16a}\\ &&
\leq 4\int_{ R^{ 2}} e^{-c{ s \over (b_n^{ -1}n)}|u|^{2} }|u|^{h/2}\,du\nn\\
&&
\leq C\({ s \over b_n^{ -1}n}\)^{ -(1+h/4)}.\nn
\eea Similarly
\begin{equation}
\int_{v+C_{ n}}|\wt{u}|^{ h/2}\,du\leq C(b_n^{
-1}n)^{(1+h/4)}.\label{18.16an}
\end{equation}

We now bound (\ref{x6.6l}) by bounding successively the integration with
respect to $u_{1},\ldots,u_{ m} $. Consider first the
$du_{ 1}$ integral, fixing $u_{ 2},\ldots,u_{ m} $. By (\ref{x6.6m}) the
$du_{1}$ integral is over the rectangle $u_{ 2}+C_{ n}$, hence the factors
involving $u_{ 1}$ can be bounded using (\ref{18.16a}). Proceeding
inductively, using (\ref{x6.6m}) when $n_{j}-n_{  j-1}>0$ and
(\ref{18.16an}) when $n_{j}=n_{  j-1}$, leads to the following bound of
(\ref{x6.6l}), and hence of (\ref{x6.6j}) on $D_{ m}( A)$:
\bea && c^{ m}\ep^{ m/2}\int_{(b_n^{ -1} n)^{ 1/2} [-\pi,\pi]^{2m}}
      \prod_{ j=1}^{ m} e^{-c|\wt{u}_{ j}|^{ 2}(n_{j}-n_{j-1})/(b_n^{-1}n) }
      |p_{ j}|^{1/2}
\,dp_{ j}\label{x6.6ja}\\ &&\leq C^{ m}\ep^{m/2}
      \prod_{ j\in A}(b_n^{ -1}n)^{ (1+ h( j)/4)}
      \prod_{ j\in A^{ c}}\({ (n_{j}-n_{j-1})
\over b_n^{ -1}n}\)^{ -(1+ h( j)/4)} .\nn
\eea Here $A^{ c}$ means the complement of $A$ in $\{ 1,\ldots,m\}$, so
that $A^{ c}$ always contains
$1$.

Note that
\begin{equation} (b_n^{ -1}n)^{ (1+ h( j)/4)/2}{1\over b_n^{-1}n}\rar 0
\mbox{ as }n\rar\ff.\label{18.46}
\end{equation} If $A=\{ i_{ 1},\ldots,i_{ k}\}$ where $i_{ 1}<\cdots<i_{
k}$ we then obtain for the sum in (\ref{18.13rx}) over
$D_{ m}( A)$, the bound
\begin{eqnarray}  &&  C^{ m} \ep^{m/4  }
\max_{h( j)}\sum_{1\leq  n_{i_{ 1}}<\cdots<n_{i_{ k}}\leq n}
      \prod_{ j\in A^{ c}}\({ (n_{ j}-n_{ j-1}) \over b_n^{ -1}n}\)^{-(1+ h(
j)/4)/2}{1\over b_n^{ -1}n}\label{18.16bb}\\ &&
\leq  C^{ m} \ep^{m/4  }
\max_{h( j)}\int_{0\leq  r_{ i_{ 1}}<\cdots<r_{i_{ k}}\leq b_n}
      \prod_{ j\in A^{ c}}(r_{j}-r_{j-1})^{-(1/2+ h( j)/8) }
\,dr_{ j}\nn\\
      &&
\leq  C^{ m} \ep^{m/4  }\max_{h( j)} {b_n^{\sum_{j\in A^{ c}}(1/2- h(
j)/8) }\over
\Ga(\sum_{ j\in A^{ c}}(1/2- h( j)/8))}
\nonumber
\end{eqnarray}

      Using this together with (\ref{18.13rx}), but with $m$ replaced by
$2m$, and the fact that $( 2m!)^{ 1/2}/m!\leq 2^{ m}$, we see that
(\ref{18.16c2})  is bounded by
\begin{eqnarray} &&
\sum_{ m=0}^{ \ff}C^{ m} \la^{ m} \ep^{m/4 }\(\sum_{A\subseteq \{
2,3,\ldots,2m\} } \max_{h( j)} {b_n^{\sum_{j\in A^{ c}}(1/2- h( j)/8)
}\over
\Ga(\sum_{ j\in A^{ c}}(1/2- h( j)/8))}\)^{ 1/2}.
\label{18.311}
\eea We have $\sum_{A\subseteq \{1, 2,3,\ldots,2m\} }1=2^{ 2m}$, and
the number of ways to choose the $\{  h( j)\}$ is bounded by the number
of ways of dividing $2m$ objects into $3$  groups, which is $ 3^{ 2m}$.
      Then noting that $\sum_{j\in A^{ c}}(1/2- h( j)/8)  $ is an integer
multiple of $1/8$ which is always less than $m$, we can bound the last
line by
\bea &&
      \sum_{ l=0}^{ \ff}\(\sum_{ m=l}^{ \ff}C^{ m}\la^{ m}\ep^{m/4 }\)
\sum_{ j=0}^{ 7}\({b_n^{l+j/8}\over \Ga(l+j/8)}\)^{ 1/2}
\label{18.311k}\\ && \leq Cb_n\sum_{ l=0}^{ \ff}\(\sum_{ m=l}^{
\ff}C^{ m}\la^{ m}\ep^{m/4 }\)
      \({b_n^{l}\over \Ga(l)}\)^{ 1/2}
\nonumber\\ &&
\leq Cb_n( 1-C\la \ep^{m/4 })^{ -1}\sum_{ l=0}^{ \ff}C^{l}\la^{
l}|\ep|^{l/4 }b_n^{l/2}
\({1\over \Ga(l)}\)^{ 1/2}\nn
\end{eqnarray} for $\ep>0$ sufficiently small.

      (\ref{18.16c2}) then follows from the fact that for any $a>0$
\begin{eqnarray}&&
\sum_{ l=0}^{ \ff}a^{l}
\({1\over \Ga(l)}\)^{ 1/2}\label{bound.k}\\ &&= \sum_{ m=0}^{
\ff}\left( a^{2m}\({1\over \Ga(2m)}\)^{ 1/2}+a^{2m+1}\({1\over
\Ga(2m+1)}\)^{ 1/2}\right)
\nn\\&&\leq   C(1+ a)\sum_{ m=0}^{ \ff} a^{2m}\({1\over
\Ga(2m)}\)^{ 1/2} \nonumber
\\&&\leq C(1+ a)e^{Ca^{ 2}}.\nonumber
\end{eqnarray}
\qed

      \begin{remark} \label{firstremark} {\rm It follows from the 
proof that in fact for $\rho>0$
      sufficiently small, for any $\la>0$
\begin{eqnarray} &&
\limsup_{ \ep\rar 0}\,\limsup_{ n\rar
\ff}\,\sup_{y}\label{18.20gen}\\ &&{ 1\over b_n}\log E\(\exp\lc
\la \Bigg|  {I_{ n}( y)-I_{n,\ep}( y)
\over \ep^{ \rho}\,b_n^{ -1}n}\Bigg|^{1/2}\rc\)=0.\nn
\end{eqnarray} }
\end{remark}

\begin{remark}\label{largedev}
{\rm Without the the restriction that $b_n=o(n)$, Theorem 
\ref{theo-md} is not true. To see this,
let $N$ be an arbitrarily large integer, let $\eps=2/N^2$,
  and let
  $X_i$ be be an i.i.d.\ sequence of random vectors
in $\Z^2$  that take the values $(N,0), (-N,0), (0,N)$, and $ (0, 
-N)$ with probability
$\eps/4$ and  $\P(X_1=(0,0))=1-\eps$. The covariance matrix of the $X_i$
will be the identity. Let $b_n=(1-\eps)n$. Then the
event that $S_i=S_0$ for all $i\leq n$ will have probability at
  least $(1-\eps)^n$, and on this event $B_n=n(n-1)/2$. This shows that
\[ \log \P(B_n-\E B_n>nb_n/2) \geq n\log(1-\eps), \]
  which would contradict (1.4).

The same example shows that the critical constant in the polymer 
model is different
than the one in \cite{BC}.
Then
\[
\E\exp\Big\{C{B_n-\E B_n\over n}\Big\}\ge
\exp\Big\{-C{\E B_n\over n}\Big\}(1-\eps)^n\exp\Big\{C{n-1\over 2}\Big\}.
\]
This shows that the critical constant is no more than
$ 2\log {1\over 1-\eps} $.  }
\end{remark}

\section{Theorem \ref{theo-ld}: Upper bound for $\E B_n-B_n$}

\medskip

{\bf  Proof of Theorem \ref{theo-ld}.}

We first prove (\ref{1.3}) for $\theta =1$:
\begin{eqnarray} -C_1&\le&
\liminf_{n\to\infty}b_n^{-1}\log\P\Big\{\E B_n -B_n
\ge (2\pi)^{-1}\det(\Gamma)^{-1/2}n\log b_n\Big\}
\nn\\ &\le&
\limsup_{n\to\infty}b_n^{-1}\log\P\Big\{\E B_n -B_n
\ge (2\pi)^{-1}\det(\Gamma)^{-1/2}n\log b_n\Big\}\le -C_2
\label{4.1}
\end{eqnarray}
      for any $\{b_n\}$ satisfying (\ref{restrict}).
$b_n^\theta$ in (\ref{4.1}).
\medskip

In this section we prove the upper bound for (\ref{4.1}). Let $t>0$ and
write
$K=[t^{-1}b_n]$. Divide $[1,n]$ into $K>1$ disjoint subintervals
$(n_0,n_1],\cdots, (n_{K-1},n_K]$, each of length $[n/K]$ or
$[n/K]+1$.  Notice that
\begin{eqnarray} &&
\E B_n- B_n\le\sum_{i=1}^K\Big[\E B\big((n_{i-1}, n_i]^2_{ <}\big)
-B\big((n_{i-1}, n_i]^2_{ <}\big)\Big]\label{4.5}\\ &&\hspace{ 2in} +\E
B_n -\sum_{i=1}^K\E B\big((n_{i-1}, n_i]^2_{ <}\big)
\nn
\end{eqnarray}

By (\ref{4.3s}),
\begin{eqnarray} &&\sum_{i=1}^K\E B\big((n_{i-1}, n_i]^2_{ <}\big)
=\sum_{i=1}^K\E B_{n_i-n_{i-1}}
\label{4.6}\\ && =\sum_{i=1}^K \Big[{1\over (2\pi)\sqrt{\det
\Gamma}}(n/K)\log(n/K) +O(n/K)\Big]\nonumber
\\ && ={1\over (2\pi)\sqrt{\det \Gamma}}n\log (n/K)+O(n)\nonumber
\end{eqnarray} With $K>1$, the error term can be taken to be
independent of
$t$ and $\{b_n\}$.
      Thus, by (\ref{4.3s}), there is constant $\log a>0$ independent of
$t$ and $\{b_n\}$ such that
\begin{eqnarray} &&\E B_n -\sum_{j=1}^K\E B\big((n_{i-1}, n_i]^2_{
<}\big)
\label{4.7}\\ && \leq {1\over (2\pi)\sqrt{\det \Gamma}}n\Big(\log
(t^{-1}b_n) +\log a\Big).\nonumber
\end{eqnarray} It is here that we use the condition that  $\E|S_{ 1}|^{
2+\de}<\ff$ for some
$\de>0$, needed for (\ref{4.3s}).

      By first using Chebyshev's inequality, then using (\ref{4.5}), (\ref{4.7})
and the independence of the $B\big((n_{i-1}, n_i]^2_{ <}\big)$, for any
$\phi >0$,
\begin{eqnarray} &&\hspace{ .5in}\P\Big\{\E B_n -B_n
\ge (2\pi)^{-1}\det(\Gamma)^{-1/2}n\log b_n\Big\}
\label{4.8}\\ && \le
\exp\Big\{-\phi b_n\log b_n\Big\}\E\exp\Big\{-2\pi\phi
\sqrt{\det \Gamma}{b_n\over n}(B_n-\E B_n)\Big\}\nonumber
\\ &&\le
\exp\Big\{\phi b_n(\log a-\log t)\Big\}
\bigg(\E\exp\Big\{-2\pi\phi
\sqrt{\det \Gamma}{b_n\over n}(B_{[n/K]}-\E B_{[n/K]})\Big\}\bigg)^K
\nonumber
\end{eqnarray} By \cite[Theorem 1.2]{R},
\begin{equation}
\sqrt{\det \Gamma}{b_n\over n}(B_{[n/K]}-\E B_{[n/K]})
\buildrel d\over\longrightarrow \gamma_t,\hskip.2in (n\to\infty)
\label{4.9}
\end{equation}
      where $\gamma_t$
      is the renormalized self-intersection local time of planar Brownian
motion
$\{ W_s\}$ up to time $t$. By Lemma \ref{lem-idbound} and the
dominated convergence theorem,
\begin{equation}
\E\exp\Big\{-2\pi\phi
\sqrt{\det \Gamma}{b_n\over n}(B_{[n/K]}-\E B_{[n/K]})\Big\}
\longrightarrow\E\exp\Big\{-2\pi\phi t\gamma_1\Big\},
\hskip.2in (n\to\infty)
\label{4.11}
\end{equation} where we used the scaling
$\gamma_t\stackrel{d}{=}t\gamma_1$.

Thus,
\begin{eqnarray} &&\limsup_{n\to\infty}b_n^{-1}\log\P\Big\{\E B_n
-B_n
\ge (2\pi)^{-1}\det(\Gamma)^{-1/2}n\log b_n\Big\}
\label{4.12}\\ &&\le \theta (\log a -\log t) +{1\over
t}\log\E\exp\Big\{-2\pi\phi t\gamma_1\Big\} \nonumber
\\ &&=\theta\log (a\phi) +{1\over t}\log\E\exp\Big\{-(\phi t)\log (\theta
t)- 2\pi(\phi t)\gamma_1\Big\} \nonumber
\end{eqnarray}
      By \cite[p. 3233]{BC}, the limit
\begin{equation} C\equiv\lim_{t\to\infty}{1\over
t}\log\E\exp\Big\{-t\log t- 2\pi t\gamma_1\Big\}
\label{4.12a}
\end{equation} exists. Hence
\begin{eqnarray} &&\limsup_{n\to\infty}b_n^{-1}\log\P\Big\{\E B_n
-B_n
\ge (2\pi)^{-1}\det(\Gamma)^{-1/2}n\log b_n\Big\}
\label{4.13}\\ &&\hspace{ 3in} \le\phi\log(a\phi) +C\phi.\nonumber
\end{eqnarray} Taking the minimizer $\phi =a^{-1}e^{-(1+C)}$ we have
\begin{eqnarray} &&\limsup_{n\to\infty}b_n^{-1}\log\P\Big\{\E B_n
-B_n
\ge (2\pi)^{-1}\det(\Gamma)^{-1/2}n\log b_n\Big\}
\label{4.14}\\ &&\hspace{ 3in}\leq -a^{-1}e^{-(1+C)}.\nonumber
\end{eqnarray} This proves the upper bound for (\ref{4.1}).

      \section{Theorem \ref{theo-ld}: Lower bound for $\E B_n-B_n$}

In this section we complete the proof of Theorem \ref{theo-ld} by proving
the lower bound for (\ref{4.1}).

Let $B(x,r)$ be the ball of radius $r$ centered at $x$. Let
$\F_k=\sigma\{ X_i: i\leq k\}$. Let us assume for simplicity that the
covariance matrix for the random walk is the identity; routine
modifications are all that are needed for the general case. We write
$\Theta$ for $(2\pi)^{-1} \hbox{ det }(\Gamma)^{-1/2}=(2\pi)^{-1}$. We
write $D( x,r)$ for the disc of radius $r$ in $\Z^{ 2}$ centered at $x$.

Let $K=[b_n]$ and $L=n/K$. Let us divide $\{1, 2, \ldots, n\}$ into
$K$ disjoint contiguous blocks, each of length strictly between $L/2$ and
$3L/2$. Denote the blocks
$J_1, \ldots, J_K$. Let  $v_i=\#(J_i)$, $w_i=\sum_{j=1}^i v_j$. Let
\begin{equation} B^{ ( i)}_{ v_{ i}}=\sum_{j,k\in J_i, j<k} \de( S_{ j},S_{
k}),
\qq A_i=\sum_{j\in J_{i-1},k\in J_{i}} \de( S_{ j},S_{ k}).\label{5.1}
\end{equation}
      Define the following sets:
\begin{eqnarray*} F_{i,1}&=&\{ S_{w_i}\in D(i\sqrt L, \sqrt L/16)\},
\label{5.2}\\ F_{i,2}&=&\{S(J_i)\subset [(i-1)\sqrt L-\sqrt L/8,i\sqrt
L+\sqrt L/8]\times  [-\sqrt L/8, \sqrt L/8]\},\nonumber\\ F_{i,3}
&=&\{B^{ ( i)}_{ v_{ i}}-\E B^{ ( i)}_{ v_{ i}}\leq
\kappa_1 L\}, \nonumber\\
      F_{i,4} &=&\{ \sum_{j\in J_{i}} 1_{D(x,r\sqrt L)}(S_j)\leq
\kappa_2 rL
\hbox{ for all } x\in D(i\sqrt L, 3\sqrt L), 1/\sqrt L<r<2\},\nonumber\\
F_{i,5} &=&\{ A_{i}<\kappa_3 L\},\nonumber
\end{eqnarray*} where
$\kappa_1, \kappa_2, \kappa_3$ are constants that will be chosen later
and do not depend on $K$ or $L$. Let
\begin{equation}C_i=F_{i,1}\cap F_{i,2}\cap F_{i,3}\cap F_{i,4}\cap
F_{i,5}\label{5.3}
\end{equation}
      and
\begin{equation}E=\cap_{i=1}^K C_i.\label{5.4}
\end{equation}

We want to show
\begin{equation}\P(C_i\mid \F_{w_{i-1}})\geq c_1>0
\label{5.5}
\end{equation}
      on the event
$C_1\cap \cdots \cap C_{i-1}$. Once we have (\ref{5.5}), then
\begin{equation}\P(\cap_{i=1}^m C_i)=\E\Big(\P(C_m\mid
\F_{w_{m-1}});
\cap_{i=1}^{m-1} C_i\Big)
\geq c_1 \P(\cap_{i=1}^{m-1} C_i),
\label{5.6}
\end{equation}
      and by induction
\begin{equation}\P(E)=\P(\cap_{i=1}^K C_i)\geq c_1^K=e^{K\log
c_1}=e^{-c_2K}.
\label{5.7}
\end{equation}

On the set $E$, we see that $S(B_i)\cap S(B_j)=\emptyset$ if $|i-j|>1$. So
we can write
\begin{equation} B_n=\sum_{k=1}^K (B^{ ( k)}_{ v_{ k}}-\E B^{ ( k)}_{ v_{
k}})+\sum_{k=1}^K \E B^{ ( k)}_{ v_{ k}}+\sum_{k=1}^K A_{k}.
\label{5.8}
\end{equation} On the event $E$, each $B^{ ( k)}_{ v_{ k}}-\E B^{ ( k)}_{
v_{ k}}$ is bounded by $\kappa_1 L$ and each $A_k$ is bounded by
$\kappa_3 L$. By (\ref{4.3}), each $\E B^{ ( k)}_{ v_{ k}}=\Theta v_k \log
v_k+O(L)=\Theta v_k \log L+O(v_k)$. Therefore
\begin{equation}B_n\leq \kappa_1KL +\Theta KL\log L+O(n)+\kappa_3KL,
\label{5.9}
\end{equation} and using (\ref{4.3}) again,
\begin{eqnarray} \E B_n- B_n &\geq& \Theta n\log n-c_3n-\Theta n\log
(n/b_n)
\label{5.9a}\\ &=&\Theta n\log b_n-c_3n\nonumber
\end{eqnarray}
      on the event $E$. We conclude that
\begin{equation}\P(\E B_n-B_n\geq \Theta n\log b_n-c_3n)\geq
e^{-c_2b_n}.
\label{5.10}
\end{equation}

We apply (\ref{5.10}) with $b_n$ replaced by $b_n'=c_4b_n$, where
$\Theta
\log c_4=c_3$. Then
\begin{equation}\Theta n\log b'_n-c_3n=\Theta n\log b_n+ \Theta n\log
c_4-c_3n=\Theta n\log b_n.
\label{5.11}
\end{equation} We then obtain
\begin{equation}\hspace{ .5in}
\P(\E B_n-B_n\geq \Theta n\log b_n)=\P(\E B_n-B_n\geq \Theta n\log
b'_n-c_3n)
\geq e^{-c_2b_n'},
\label{5.12}
\end{equation} which would complete the proof of the lower bound for
(\ref{4.1}), hence of Theorem \ref{theo-ld}.

So we need to prove (\ref{5.5}). By scaling and the support theorem for
Brownian motion (see \cite[Theorem I.6.6]{B}), if $W_t$ is a planar
Brownian motion and $|x|\leq \sqrt L/16$, then
\begin{eqnarray} &&\hspace{ .2in}\P^x\(W_{v_i}\in D(\sqrt L, \sqrt L/16)
\mbox{ and}\right.
\label{5.13}\\ &&\hspace{ .6in}\left.\{W_s; 0\leq s\leq {v_i}\}
\subset [-\sqrt L/8, 9\sqrt L/8]\times [-\sqrt L/8, \sqrt L/8]\) >  c_5,
\nonumber
\end{eqnarray} where $c_5$ does not depend on $L$. Using Donsker's
invariance principle for random walks with finite second moments
together with the Markov property,
\begin{equation}\P(F_{i,1}\cap F_{i,2}\mid F_{w_{i-1}})> c_6.
\label{5.14}
\end{equation} By Lemma \ref{lem-idbound}, for $L/2\leq \ell\leq 3L/2$
\begin{equation}\P(B_\ell-\E B_\ell> \kappa_1 L)\leq c_6/2\label{5.14a}
\end{equation} if we choose $\kappa_1$ large enough. Again using the
Markov property,
\begin{equation}\P(F_{i,1}\cap F_{i,2}\cap F_{i,3}\mid F_{w_{i-1}})>
c_6/2.\label{5.15}
\end{equation}

Now let us look at $F_{i,4}$. By \cite[ p. 75]{S}, $\P(S_j=y)\leq c_7/j$ with
$c_7$ independent of $y\in \Z^{ 2}$ so that
\begin{equation}
\P(S_j\in D(x, r\sqrt L))=\sum_{ y\in D(x, r\sqrt L)}\P(S_j=y)\leq {c_8r^{
2}L
\over j}.\label{5.15s}
\end{equation} Therefore
\begin{eqnarray}\hspace{ .2in}\E\sum_{j\in J_1} 1_{D(x,r\sqrt L)}(S_j)
&\leq&
\sum_{j=1}^{[2L]} \P(S_j\in D(x, r\sqrt L))
\label{}\\ &\leq& r^2L+\sum_{j=r^2L}^{[2L]}
\frac{c_{9}r^2L}{j}\nn\\ &\leq& r^2L+c_{10}Lr^2 \log(1/r)\le
c_{11}Lr^2 \log(1/r)\nonumber
\end{eqnarray}
      if
$1/\sqrt L\leq r\leq 2$. Let $C_m=\sum_{j< m} 1_{D(x,r\sqrt L)}(S_j)$ for
$m\leq {[2L]}+1$ and let $C_m=C_{{[2L]}+1}$ for
$m>L$. By the Markov property and independence,
\begin{eqnarray} &&\E[C_\infty-C_m\mid \F_m]\leq
1+\E[C_\infty-C_{m+1}\mid
\F_m]
\label{5.17}\\ &&\hspace{ 1.3in}\leq 1+ \E^{S_m}C_\infty\leq
c_{12}Lr^2\log (1/r).
\nonumber
\end{eqnarray} By \cite[Theorem I.6.11]{B}, we have
\begin{equation}\E \exp\Big(c_{13}
\frac{C_{{[2L]}+1}}{c_{12}Lr^2\log (1/r)}\Big)\leq c_{14}\label{5.18}
\end{equation} with $c_{13}, c_{14}$ independent of $L$ or
$r$. We conclude
\begin{equation}\P\Big(\sum_{j\in J_1} 1_{D(x,r\sqrt
L)}(S_j)>c_{15}Lr^2\log(1/r)\Big)\leq c_{16}e
^{-c_{17}c_{15}}.\label{5.19}
\end{equation} Suppose $2^{-s}\leq r<2^{-s+1}$ for some $s\geq 0$. If
$x\in D(0,3\sqrt L)$, then each point in the disc $D(x,r\sqrt L)$ will be
contained in
$D(x_i,  2^{-s+3}\sqrt L)$ for some $x_i$, where each coordinate of
$x_i$ is an integer multiple of $2^{-s-2}\sqrt L$. There are at most
$c_{18}2^{2s}$ such balls, and $Lr^2\log(1/r)\leq c_{ 19}2^{ s/2}Lr$, so
\begin{equation}\hspace{ .3in}\P\Big(\sup_{x\in D(0,3\sqrt L), 2^{-s}\leq
r<2^{-s+1}}
\sum_{j\in J_1} 1_{D(x,r\sqrt L)}(S_j)> c_{20}rL\Big)\leq c_{21}2^{2s}
e^{-c_{22}c_{20} 2^{s/2}}.
\label{5.20}
\end{equation} If we now sum over positive integers $s$ and take
$\kappa_2=c_{20}$ large enough, we see that
\begin{equation}\P(F_{1,4})\leq c_6/4.
\label{5.21}
\end{equation} By the Markov property, we then obtain
\begin{equation}\P(F_{i,1}\cap F_{i,2}\cap F_{i,3} \cap F_{i,4}\mid
F_{w_{i-1}})> c_6/4.\label{5.22}
\end{equation}

Finally, we examine $F_{i,5}$. We will show
\begin{equation}\P(F_{i,5}\mid \F_{w_{i-1}})\leq c_6/8
\label{5.23}
\end{equation} on the set
$\cap_{j=1}^{i-1} C_j$ if we take $\kappa_3$ large enough. By the Markov
property, it suffices  to show
\begin{equation}\P\Big( \sum_{j=1}^{[2L]} 1_{(S_j\in G)} \geq
\kappa_3L\Big)\leq c_6/8
\label{5.24}
\end{equation}
      whenever $G\in \Z^{ 2}$ is a fixed nonrandom set consisting of
${[2L]}$ points satisfying the property that
\begin{equation}\hspace{ .3in}\#(G\cap D(x,r\sqrt L))\leq
\kappa_2rL, \qq x\in D(0,3\sqrt L),
\quad 1/\sqrt L\leq r\leq 2.
\label{5.25}
\end{equation}
      We compute the expectation of
\begin{equation}\sum_{j=1}^{[2L]} 1_{(S_j\in G\cap (D(0, 2^{-k}\sqrt
L)\setminus D(0, 2^{-k+1}\sqrt L)))}.
\label{5.26}
\end{equation} When $j\leq 2^{-2k}L$, then the fact that the random walk
has finite second moments implies that
      the probability that $|S_j|$ exceeds $2^{-k+1}\sqrt L$ is bounded by
$c_{23}j/(2^{-2k+2}L)$. When $j> 2^{-2k} L$, we use \cite[p. 75]{S}, and
obtain
\begin{equation}\P(S_j\in G\cap (D(0, 2^{-k}\sqrt L)\leq c_{24}
\frac{\kappa_2 2^{-k}L}{j}.
\label{5.27}
\end{equation}
      So
\begin{eqnarray} &&\E \sum_{j=1}^{[2L]} 1_G(S_j)
\label{5.28}\\ && \leq
\sum_{k}\sum_{{[2L]}\geq j>2^{-2k}L} c_{24} \frac{\kappa_2 2^{-k}L}{j}
+\sum_k\sum_{j\leq 2^{-2k}L} c_{23}
\frac{j}{2^{-2k+2}L}\nonumber\\ &&\leq \sum_k(c_{25} \kappa_2
k2^{-k}L+c_{26}2^{-2k}L)
\leq c_{27}L.\nn
\end{eqnarray}
      So if take $\kappa_3$ large enough, we obtain (\ref{5.24}).

This completes the proof of (\ref{5.5}), hence of Theorem
\ref{theo-ld}.
\qed

\section{ Laws of the iterated logarithm}

\subsection{Proof of the LIL for $B_n-\E B_n$}

First, let $S_j, S'_j$ be two independent copies of our random walk. Let
\begin{equation}
\ell(n,x)=\sum_{i=1}^n \de( S_i\,,x),\hspace{ .5in}\ell'(n,x)
=\sum_{i=1}^n
\de( S'_i\,,x)\label{6.1a}
\end{equation} and note that
\begin{equation}I_{k,n }=\sum_{i=1}^k \sum_{j=1}^n \de(S_i,S'_j)
=\sum_{x\in \Z^2}\ell(k,x)\ell'(n,x).\label{6.1}
\end{equation}

\bl \label{lem-lilexp}There exist constants $c_1, c_2$ such that
\begin{equation}\P(I_{k,n }> \lambda \sqrt{kn})\leq c_1
e^{-c_2\lam}.\label{6.2}
\end{equation}
\el

\proof  Clearly
\begin{equation} (I_{k,n })^m =\sum_{x_1\in \Z^2}\cdots
\sum_{x_m\in \Z^2}
\Big(\prod_{i=1}^m \ell(k, x_i)\Big) \Big(\prod_{i=1}^m \ell'(n, x_i)\Big)
\label{6.3}
\end{equation}
      Using the independence of $S$ and $S'$,
\begin{equation}
\E\( (I_{k,n })^m\)=\sum_{x_1\in \Z^2}\cdots \sum_{x_m\in \Z^2}
\E\Big(\prod_{i=1}^m \ell(k, x_i)\Big)\E \Big(\prod_{i=1}^m \ell'(n,
x_i)\Big).
\label{6.4}
\end{equation} By Cauchy-Schwarz, this is less than
\begin{eqnarray} &&  \Big[\sum_{x_1\in \Z^2}\cdots \sum_{x_m\in
\Z^2}
\Big(\E\Big(\prod_{i=1}^m \ell(k,
x_i)\Big)\Big)^2\Big]^{1/2}\label{6.5}\\ &&\hspace{ 1in}
      \Big[\sum_{x_1\in \Z^2}\cdots \sum_{x_m\in \Z^2}
\Big(\E\Big(\prod_{i=1}^m \ell'(n, x_i)\Big)\Big)^2\Big]^{1/2}
\nn\\ &&\hspace{2in}=: J_1^{1/2} J_2^{1/2}. \nonumber
\end{eqnarray}
      We can rewrite
\begin{equation}J_1=\sum_{x_1\in \Z^2}\cdots \sum_{x_m\in \Z^2}
\E\Big(\prod_{i=1}^m \ell(k, x_i)\Big)\E \Big(\prod_{i=1}^m \ell'(k,
x_i)\Big)=\E\( (I_{k})^m\),
\label{6.6}
\end{equation} and  similarly
$J_2=\E\( (I_{n})^m\)$.

Therefore,
\begin{eqnarray} &&\E \exp(aI_{k,n }/\sqrt{kn})\label{6.7}\\
&=&\sum_{m=0}^\infty
\frac{a^m}{k^{m/2}n^{m/2} m!} \E\( (I_{k,n})^m\)
\nn\\ &\leq&
\sum_m\frac{a^m }{k^{m/2}n^{m/2} m!} (\E\( (I_{k})^m\))^{1/2} (\E\(
(I_{n})^m\))^{1/2} \nonumber\\ &\leq&\Big(\sum \frac{a^m }{m!} \E
\Big(\frac{ I_{k}}{k}\Big)^m\Big)^{1/2}
      \Big(\sum \frac{a^m }{m!} \E \Big(\frac{
I_{n}}{n}\Big)^m\Big)^{1/2}\nn\\ &\leq&\Big(\E e^{a I_{k}/k}\Big)^{1/2}
\Big(\E e^{a   I_{n}/n}\Big)^{1/2}.\nn
\end{eqnarray} By Lemma \ref{lem-ubound} this can be  bounded
independently of $k$ and $n$ if $a$ is taken small, and our result follows.
\qed

We are now ready to  prove the upper bound for the LIL for $B_n-\E B_n$.
Write
$\Xi$ for $\sqrt{\det \Gamma}\,\,\kappa(2,2)^{-4}$. Recall that for any
integrable random variable $Z$ we let $\ol Z$ denote $Z-\E Z$. Let
$\eps>0$ and let $q>1$ be chosen later.  Our first goal is to get an upper
bound on
$$\P(\max_{n/2\leq k\leq n} \ol B_k>(1+\eps)\Xi^{-1} n\log\log n).$$ Let
$m_0=2^N$, where $N$ will be chosen later to depend only on
$\eps$. Let $\sA_0$ be the integers of the form $n-km_0$ that are
contained in $ \{n/4, \ldots, n\}$.  For each $i$ let $\sA_i$ be the set of
integers of the form $n-km_02^{-i}$ that are contained in
$\{n/4, \ldots, n\}$. Given an integer $k$, let $k_j$ be the largest element
of
$\sA_j$ that is less than or equal to $k$. For any $k\in
\{n/2, \ldots, n\}$, we can write
\begin{equation}\ol B_k=\ol B_{k_0}+(\ol B_{k_1}-\ol B_{k_0})+\cdots+
(\ol B_{k_N}-\ol B_{k_{N-1}}).\label{6.8}
\end{equation}

\ni If $\ol B_k\geq (1+\eps)\Xi^{-1}n\log\log n$ for some $n/2\leq k\leq
n$, then either \hfill\break
\qq (a) $\ol B_{k_0}\geq (1+\frac{\eps}{2})\Xi^{-1} n\log\log n$ for some
$k_0\in \sA_0$; or else\hfill\break \qq (b) for some $i\geq 1$ and some
pair of consecutive elements $k_i, k'_i
\in \sA_i$, we have
\begin{equation} \ol B_{k'_i}-\ol B_{k_i}\geq \tfrac{\eps}{40i^2}
n\log\log n.
\label{6.9}
\end{equation}

For each $k_0$,  using Theorem \ref{theo-md} and the fact that
$k_0\geq n/4$, the probability in (a) is bounded by
\begin{equation}
\exp(-(1+\tfrac{\eps}{4})\log\log k_0)\leq c_1 (\log
n)^{-(1+\frac{\eps}{4})}.
\label{6.10}
\end{equation}
       There are at most $n/m_0$ elements of
$\sA_0$, so the probability in (a) is bounded by
\begin{equation}
\frac{n}{m_0} \frac{c_1}{(\log n)^{1+\frac{\eps}{4}}}.\label{6.11}
\end{equation}

Now let us examine the probability in (b). Fix $i$ for the moment. Any two
consecutive elements of $\sA_i$ are $2^{-i}m_0$ apart. Recalling the
notation (\ref{2.9}) we can write
\begin{equation}\ol B_k-\ol B_j=\ol B([j+1,k]^{ 2}_{ <})+\ol B([1,j]\times
[j+1,k]),
\label{6.12}
\end{equation}
      So
\begin{eqnarray} \P(\ol B_k-\ol B_j\geq  \tfrac{\eps}{40i^2} n\log
\log n)&\leq &\P(\ol B([j+1,k]^{ 2}_{ <})\geq
\tfrac{\eps}{80i^2} n\log \log n)\nn
\\ & &\hspace{ -.25in} +\P\(B([1,j]\times [j+1,k])\geq\tfrac{\eps}{80i^2}
n\log \log n\).\label{6.14}
\end{eqnarray} We bound the first term on the right by Theorem
\ref{theo-md}, and get the bound
\begin{equation}\exp\Big(-\frac{\eps \Xi}{80i^2} \frac{n\log\log
n}{2^{-i}m_0}\Big)
\leq \exp\Big(-\frac{\eps \Xi}{80i^2} 2^i(n/m_0)\log\log n\Big)
\label{6.15}
\end{equation} if $j$ and $k$ are consecutive elements of $\sA_i$. Note
that
$B([1,j]\times [j+1,k])$ is equal in law to
$I_{j-1,k-j}$. Using Lemma
\ref{lem-lilexp}, we bound the second term on the right hand side of
(\ref{6.14}) by
\begin{eqnarray} &&c_1\exp\Big(-c_2\frac{\eps}{80i^2}\frac{n\log\log
n}{\sqrt{2^{-i}m_0}\sqrt j}\Big)\nn\\ &&
\leq c_1\exp\Big(-c_2\frac{\eps}{80i^2} 2^{i/2} (n/m_0)^{1/2}\log
\log n\big).
\label{6.16}
\end{eqnarray} The number of pairs of consecutive elements of
$\sA_i$ is less than $2^{i+1}(n/m_0)$. So if we add (\ref{6.15})  and
(\ref{6.16}) and multiply by the number of pairs, the probability of (b)
occurring for a fixed
$i$ is bounded by
\begin{equation}c_3 \frac{n}{m_0} 2^i \exp\Big(-c_4 2^{i/2}
(n/m_0)^{1/2}\log\log n/(80i^2)\Big).
\label{6.16a}
\end{equation}
      If we now sum over $i\geq 1$, we bound the probability in (b) by
\begin{equation}c_5 \frac{n}{m_0}\exp\Big(-c_6 (n/m_0)^{1/2}\log\log
n\Big).
\label{6.17}
\end{equation} We now choose $m_0$ to be the largest power of
$2$ so that $c_6 (n/m_0)^{1/2}>2$; recall  $n$ is big.

Let us use this value of $m_0$ and combine (\ref{6.11}) and (\ref{6.17}).
Let
$n_\ell=q^\ell$ and
\begin{equation}C_\ell=\{ \max_{n_{\ell-1}\leq k\leq n_\ell} \ol B_k\geq
(1+\eps)\Xi^{-1}n_{\ell}
\log\log n_{\ell}\}.
\label{6.18}
\end{equation}
      By our estimates, $\P(C_\ell)$ is summable, so for $\ell$ large, by
Borel-Cantelli we have
\begin{equation}\max_{n_{\ell-1}\leq k\leq n_\ell} \ol B_k\leq
(1+\eps)\Xi^{-1} n_\ell \log\log n_\ell.
\label{6.19}
\end{equation} By taking $q$ sufficiently close to 1, this implies that for
$k$ large we have $\ol B_k\leq (1+2\eps)\Xi^{-1} k\log\log k$. Since
$\eps$ is arbitrary, we have our upper bound.
\qed

The lower bound for the first LIL is easier. Let $\delta>0$ be small and let
$n_\ell=[e^{\ell^{1+\delta}}]$. Let
\begin{equation}D_\ell=\{ \ol B([n_{\ell-1}+1,n_\ell]^{ 2}_{ <})
\geq (1-\eps)\Xi^{-1} n_\ell \log\log n_\ell\}.\label{6.20}
\end{equation}
      Using Theorem 1.1, and the fact that $n_\ell/(n_\ell -n_{\ell-1})\sim
ce^{ -c'l^{
\de}}$ we see that
$\sum_\ell \P(D_\ell)=\infty$. The $D_\ell$ are independent, so  by
Borel-Cantelli
\begin{equation}\ol B([n_{\ell-1}+1,n_\ell]^{ 2}_{ <})\geq
(1-\eps)\Xi^{-1} n_\ell
\log\log n_\ell
\label{6.21}
\end{equation} infinitely often with probability one.   Note that as in
(\ref{6.12})  we can write
\begin{equation}\ol B_{n_\ell}=\ol B([n_{\ell-1}+1,n_\ell]^{ 2}_{ <})+\ol
B_{n_{\ell-1}} +\ol B([1,n_{\ell-1}]\times [n_{\ell-1}+1,n_\ell]).
\label{6.22}
\end{equation} By the upper bound,
\[\limsup_{\ell\to \infty} \frac{\ol B_{n_{\ell-1}}}{n_{\ell-1}\log\log
n_{\ell-1}}\leq \Xi^{-1}
\] almost surely, which implies
\begin{equation}\limsup_{\ell\to \infty} \frac{\ol B_{n_{\ell-1}}}{n_\ell
\log\log n_\ell}=0.
\label{6.23}
\end{equation} Since $B([1,n_{\ell-1}]\times [n_{\ell-1}+1,n_\ell])\geq 0$
and by (\ref{2.0y})
\begin{equation}\hspace{ .5in}\E B([1,n_{\ell-1}]\times
[n_{\ell-1}+1,n_\ell])\leq c_1\sqrt{n_{\ell-1}}\sqrt{n_\ell-n_{\ell-1}}
=o(n_\ell\log\log n_\ell),
\label{6.24}
\end{equation} using (\ref{6.21})-(\ref{6.24}) yields the lower bound.
\qed

\subsection{ LIL for $\E B_n-B_n$}

Let $\Delta=2\pi\sqrt{\det \Gamma}$. Let us write $J_n=\E B_n-B_n$.

First we do the upper bound. Let $m_0$, $\sA_i$, and $k_j$ be as in the
previous subsection. We write, for $n/2\leq k\leq n$,
\begin{equation}J_k=J_{k_0}+(J_{k_1}-J_{k_0})+\cdots+(J_{k_N}-J_{k_{N-1}}).
\label{7.1}
\end{equation}
\ni If $\max_{n/2\leq k\leq n} J_k\geq (1+\eps)\Delta^{-1} n\log\log\log
n$, then either
\hfill\break \qq (a) $J_{k_0}\geq (1+\frac{\eps}{2}) \Delta^{-1} n
\log\log\log n$ for some
$k_0\in \sA_0$, or else \hfill \break \qq (b) for some $i\geq 1$ and
$k_i, k'_i$ consecutive elements of $\sS_i$ we have
\begin{equation}J_{k'_i}-J_{k_i}\geq \frac{\eps}{40i^2} n\log\log\log n.
\label{7.2}
\end{equation}

There are at most $n/m_0$ elements of $\sA_0$. Using Theorem
\ref{theo-ld}, the probability of (a) is bounded by
\begin{equation}c_1\frac{n}{m_0}e^{-(1+\frac{\eps}{4})\log \log n}.
\label{7.3}
\end{equation}

To estimate the probability in (b), suppose $j$ and $k$ are consecutive
elements of $\sA_i$. There are at most
$2^{i+1}(n/m_0)$ such pairs. We have
\begin{eqnarray} J_k-J_j&=&-\ol B([j+1,k]^{ 2}_{ <})-
\ol B([1,j]\times [j+1,k])\label{7.4}\\ &&\leq -\ol B([j+1,k]^{ 2}_{ <})+\E
B([1,j]\times [j+1,k])
\nn\\ &&  \leq  -\ol B([j+1,k]^{ 2}_{ <})+c_2\sqrt j\sqrt{k-j},\nonumber
\end{eqnarray}
      as in the previous subsection. Provided $n$ is large enough,
$c_2\sqrt j\sqrt{k-j}=c_2\sqrt j\sqrt{2^{-i}m_0}$ will be less than
$\frac{\eps}{80i^2}n\log\log\log n$ for all $i$. So in order for
$J_k-J_j$ to be larger than $\frac{\eps}{40i^2}n\log\log\log n$, we must
have
$-\ol B([j+1,k]^{ 2}_{ <})$ larger than
$\frac{\eps}{80i^2}n\log\log\log n$. We use Theorem 1.2 to bound this.
Then  multiplying by the number of pairs and summing over $i$, the
probability is (b) is bounded by
\begin{eqnarray} &&\sum_{i=1}^\infty 2^{i+1}\frac{n}{m_0}
e^{-\frac{\eps}{80i^2} \frac{n}{2^{-i}m_0}
\log\log n}\leq c_3\frac{n}{m_0} e^{-c_4 (n/m_0)\log\log n}.\label{7.5}
\end{eqnarray} We choose $m_0$ to be the largest possible power of 2
such that $c_4(n/m_0)>2$.

Combining (\ref{7.3}) and (\ref{7.5}), we see that if we set $q>1$ close to
1,
$n_\ell=[q^\ell]$, and \begin{equation}E_\ell=\{\max_{n_\ell/2\leq
k\leq n_\ell} J_k\geq (1+\eps) \Delta^{-1}n_\ell \log\log\log n_\ell\},
\label{7.6}
\end{equation} then
$\sum_\ell \P(E_\ell)$ is finite. So by Borel-Cantelli, the event
$E_\ell$ happens for a last time, almost surely. Exactly as in the previous
subsection, taking $q$ close enough to 1 and using the fact that $\eps$ is
arbitrary leads to the  upper bound.
      \qed

The proof of the lower bound
      is fairly similar to the previous subsection. Let $n_\ell=
[e^{\ell^{1+\delta}}]$. Theorem \ref{theo-ld} and Borel-Cantelli tell us that
$F_\ell$ will happen infinitely often, where
\begin{equation}F_\ell=\{ -\ol B([n_{\ell-1}+1,n_\ell]^{ 2}_{ <})\geq
(1-\eps)\Delta^{-1} n_\ell
\log\log\log n_\ell\}.
\label{7.7}
\end{equation}
      We have
\begin{equation}J_{n_\ell}\geq -\ol B([n_{\ell-1}+1,n_\ell]^{ 2}_{
<})+J_{n_{\ell-1}}-A(1, n_{\ell-1}; n_{\ell-1}, n_\ell).
\label{7.8}
\end{equation}
      By the upper bound,
\begin{equation} J_{n_{\ell-1}}=O(n_{\ell-1}\log \log\log
n_{\ell-1})=o(n_\ell\log\log\log n_\ell).
\label{7.9}
\end{equation}
      By Lemma \ref{lem-lilexp},
\begin{equation}\P(B([1,n_{\ell-1}]\times [n_{\ell-1}+1,n_\ell])\geq
\eps n_\ell\log\log\log n_\ell)
\leq c_1 \exp\Big(- c_2\frac{\eps n_\ell\log\log\log n_\ell}{
\sqrt{n_{\ell-1}}\sqrt{n_\ell-n_{\ell-1}}}\Big).
\label{7.10}
\end{equation}
       This is summable in
$\ell$, so
\begin{equation}\limsup_{\ell\to \infty} \frac{B([1,n_{\ell-1}]\times
[n_{\ell-1}+1,n_\ell])}{ n_\ell\log\log\log n_\ell}
\leq \eps
\label{7.11}
\end{equation} almost surely. This is true for every $\eps$, so the limsup
is 0. Combining this with (\ref{7.9}) and substituting in (\ref{7.8})
completes the proof.
\qed

\def\noopsort#1{} \def\printfirst#1#2{#1} \def\singleletter#1{#1}
       \def\switchargs#1#2{#2#1} \def\bibsameauth{\leavevmode\vrule
height .1ex
       depth 0pt width 2.3em\relax\,}
\makeatletter \renewcommand{\@biblabel}[1]{\hfill#1.}\makeatother
\newcommand{\bysame}{\leavevmode\hbox to3em{\hrulefill}\,}

\end{document}